\documentclass[10pt]{amsart}
\usepackage[T1]{fontenc}
\usepackage{geometry}
\usepackage[latin1] {inputenc}
\usepackage{amsmath}
\usepackage{amsfonts, amssymb, textcomp}
\usepackage[colorlinks=flase, linkcolor=red,urlcolor=green, citecolor=blue]{hyperref}
\usepackage{subeqnarray}

\usepackage{latexsym}
\usepackage{fancyhdr}
\usepackage{longtable}
\usepackage{amsmath, amssymb}
\usepackage{graphicx}
\usepackage{enumerate}

\numberwithin{equation}{section}

\theoremstyle{plain}
\newtheorem{proposition}{Proposition}[section]
\newtheorem{theorem}[proposition]{Theorem}
\newtheorem{lemma}[proposition]{Lemma}
\newtheorem{corollary}[proposition]{Corollary}
\newtheorem{definition}[proposition]{Definition}

\newtheorem{remark}[proposition]{Remark}

\newcommand{\RR}{\mathbb{R}}
\newcommand{\CC}{\mathbb{C}}
\newcommand{\NN}{\mathbb{N}}

\newcommand{\QQ}{\mathbb{Q}}

\let\on=\operatorname

\newenvironment{demo}[1]{\par\smallskip\noindent{\bf #1.}}{\par\smallskip}

\newsavebox{\fmbox}
\newenvironment{fmpage}[1]
 {\begin{lrbox}{\fmbox}\begin{minipage}{#1}}
 {\end{minipage}\end{lrbox}\fbox{\usebox{\fmbox}}}

\title[On the construction of large algebras]
{On the construction of large algebras not contained in the image
of the Borel map}

\author[C.~Esser]{C\'eline Esser}

\author[G.~Schindl]{Gerhard Schindl}

\address{C.~Esser: Universit\'e de Li\`ege, D\'epartement de Math\'ematique, Quartier Polytech 1, All\'ee de la D\'ecouverte 12, B\^atiment B37, B-4000 Li\`ege, Belgique}
\email{celine.esser@uliege.be}

\address{G.~Schindl: Fakult\"at f\"ur Mathematik, Universit\"at Wien, Oskar-Morgenstern-Platz~1, A-1090 Wien, Austria.}
\email{gerhard.schindl@univie.ac.at}

\begin{document}

\begin{abstract}
The Borel map $j^{\infty}$ takes germs at $0$ of smooth functions to
the sequence of iterated partial derivatives at $0$. It is well known
that the restriction of $j^{\infty}$ to the germs of quasianalytic
ultradifferentiable classes which are strictly containing the real
analytic functions can never be onto the corresponding sequence
space. In a recent paper the authors have studied the size of the
image of $j^{\infty}$ by using different approaches and worked in the
general setting of quasianalytic ultradifferentiable classes defined
by weight matrices. The aim of this paper is to show that the image of
$j^{\infty}$ is also small with respect to the notion of algebrability
and we treat both the Cauchy product (convolution) and the pointwise
product. In particular, a deep study of the stability of the
considered spaces under the pointwise product is developed.
\end{abstract}

\thanks{G. Schindl is supported by FWF-Project J~3948-N35.}
\date{\today}

\maketitle

\section{Introduction}\label{Introduction}

Classes of ultradifferentiable functions on an open subset $U\subseteq
\RR$ are classically defined by imposing growth restrictions on
their derivatives. In the case these restrictions are controlled by a
weight sequence $M=(M_{j})_{j \in \NN}$, given a sequence
$\mathbf{a}=(a_{j})_{j \in \NN}$ of complex numbers, many authors have investigated
under which conditions on $M$ and $\mathbf{a}$ there exists a function $f$ in
the class associated to $M$ satisfying $f^{(j)}(0) = a_{j}$ for every
$j \in \NN$, see \cite{Carleman23,petzsche,thilliez}. This coincides with the study of the
surjectivity of the Borel map $f \mapsto \big(f^{(j)}(0)\big)_{j \in
  \NN}$ in the corresponding spaces. Following
the work of \cite{BraunMeiseTaylor90}, it is also very classical to consider growth
restrictions defined by using weight functions $\omega$. In this
situation, the study of the surjectivity of the Borel map has been
proposed in \cite{BonetMeiseTaylorSurjectivity,BonetMeise}. More recently, new classes of ultradifferentiable
functions have been introduced in order to obtain a general framework
that covers both previous situtations, but also different ones, see \cite{compositionpaper} and \cite{dissertation}. These
classes are based on weight matrices $\mathcal{M}$ and the study of
the surjectivity of the Borel map in this context has been carried
out in \cite{borelmappingquasianalytic}. In any situation, it appears that if the considered
class is quasianalytic, which means that on this class the Borel map
is injective, and if it contains strictly the analytic functions, then
the Borel map is never surjective onto the corresponding weighted sequence space. In this context, the authors have
studied in the recent paper \cite{Borelmapgenericity} the question of knowing how
far is the Borel map from being surjective. More precisely, they
obtained that the image of the Borel map is ``small'' in the
corresponding sequence space, where the notion of smallness is defined
using different approaches: the notion of residual
sets based on Baire categories, the notion of prevalence, and the
notion of lineability. This paper aims
at obtaining the corresponging result in the algebraic sense, using
the notion of {\itshape algebrability.} While the concept of lineability consists
in proving the existence of large linear subspace satisfying a
particular property, one could search for other structure, such as
algebra, see \cite{BPS} and \cite{ABPS} and the references therein.

\begin{definition}
Let $\mathcal{A}$ be an algebra and $\kappa$ be a cardinal number. A subset $\mathcal{B}\subseteq
\mathcal{A}$ is \emph{$\kappa$-algebrable} if there is a
$\kappa$-generated subalgebra
$\mathcal{C}\subseteq \mathcal{B}\cup \{0\}$.
\end{definition}

The results of \cite{Borelmapgenericity} will be extended in two ways: first, we will
consider that the multiplicative structure on the weighted formal power series space is
given by the Cauchy (or convolution) product, which corresponds to the natural pointwise product of functions. This will be the core of
Section \ref{Algebrconvolution}. In this context, it
seems to be more natural to consider  weighted formal power series
spaces instead of sequences spaces, see Remark
\ref{onetoonepowerseries} for some explanations: this will be done in
this paper.
 In Section \ref{Algebrpointwise}, we
will work under the assumption that the multiplication is the
pointwise product. In particular, a deep study of the stability of the
image and the corresponding power series space under the pointwise product
is proposed in Section \ref{Algebrpointwise} for weight sequences
and weight matrices, and in  Section
\ref{stabilitypointwiseweightfunction} for weight functions.
 We will see that, contrary to what happens in the
case of the Cauchy product, under our assumptions, this product does
not make sense in the case of a weight sequence, or a weight
function. However, we will construct in Section \ref{Algebrpointwise}
an example of a weight matrix which gives the stability of the
corresponding space under the pointwise product which underlines the different behavior of classes defined by general weight matrices.

\medskip
Let us mention that  Section \ref{weightsequencesandgerms} is
dedicated to remind the reader the classes associated to weight sequences.  The presentation of this work and the standard assumptions on the
weight structures are similar to the ones considered in
\cite{borelmappingquasianalytic} and
\cite{Borelmapgenericity}. Moreover, throughout this paper, we write
$\NN=\{0,1,\dots\}$, $\mathcal{E}(U)$ and $\mathcal{C}^{\omega}(U)$
shall denote respectively the class of all $\CC$-valued smooth
functions and the class of all real analytic functions defined on
non-empty open $U\subseteq\RR$. For reasons of convenience we will
write $\mathcal{E}_{[M]}$ if either $\mathcal{E}_{\{M\}}$ or
$\mathcal{E}_{(M)}$ is considered, but not mixing the cases if
statements involve more than one $\mathcal{E}_{[M]}$ symbol. We use
similar notations for the classes of weighted formal power series
$\mathcal{F}_{[M]}$ and for classes defined by weight functions $\omega$ and weight matrices $\mathcal{M}$ as well. Finally, the cardinal $\mathfrak{c}$ will denote the
continuum. 

\section{Weight sequences and germs of ultradifferentiable functions}\label{weightsequencesandgerms}

\begin{definition}
Let $M=(M_j)_{j \in \NN}\in\RR_{>0}^{\NN}$ be an arbitrary sequence of
positive real numbers. Let 
$U\subseteq\RR$ be non-empty and open. The $M$-\emph{ultradifferentiable Roumieu type class} is defined by
\begin{equation*}\label{roumieu}
\mathcal{E}_{\{M\}}(U):=\{f\in\mathcal{E}(U):\,\forall\,K\subseteq U\,\text{compact}\,\, \exists\,h>0,\,\|f\|^M_{K,h}<+\infty\},
\end{equation*}
and the $M$-\emph{ultradifferentiable Beurling type class} by
\begin{equation*}\label{beurling}
\mathcal{E}_{(M)}(U):=\{f\in\mathcal{E}(U):\,\forall\,K\subseteq U\,\text{compact}\,\, \forall\,h>0,\,\|f\|^M_{K,h}<+\infty\},
\end{equation*}
where
\begin{equation*}\label{semi-norm-1}
\|f\|^M_{K,h}:=\sup_{j\in\NN,x\in K}\frac{|f^{(j)}(x)|}{h^j M_{j}}.
\end{equation*}
\end{definition}

Moreover we will write $m=(m_j)_{j \in \NN}$ for $m_j:=\frac{M_j}{j!}$.

For any compact set $K$ with smooth boundary $\mathcal{E}_{M,h}(K):=\{f\in\mathcal{E}(K): \|f\|^M_{K,h}<+\infty\}$ is a Banach space.
The Roumieu type class is endowed with the projective topology with respect to all $K\subseteq U$ compact and the inductive topology with respect to $h\in\NN_{>0}$. Similarly  the Beurling type class is endowed with the projective topology with respect to $K\subseteq U$ compact and with respect to $1/h$, $h\in\NN_{>0}$. Hence $\mathcal{E}_{(M)}(U)$ is a {\itshape Fr\'echet space} and $\underset{h>0}{\varinjlim}\,\mathcal{E}_{M,h}(K)=\underset{n\in\NN_{>0}}{\varinjlim}\,\,\mathcal{E}_{M,n}(K)$ is a {\itshape Silva space}, i.e. a countable inductive limit of Banach spaces with compact connecting mappings, see \cite[Proposition 2.2]{Komatsu73}.\vspace{6pt}

Note that the special case $M=(j!)_{j\in \NN}$ yields
$\mathcal{E}_{\{M\}}(U)=\mathcal{C}^{\omega}(U)$ the space of real
analytic functions on $U$, whereas $\mathcal{E}_{(M)}(U)$ consists of the restrictions of all entire functions provided that $U$ is connected.

\begin{definition}
The spaces of \emph{germs at $0\in\RR$} of the $M$-ultradifferentiable functions of Roumieu and Beurling types are defined respectively by
\begin{equation*}\label{roumieugerm}
\mathcal{E}_{\{M\}}^{0}:=\underset{k\in\NN_{>0}}{\varinjlim}\mathcal{E}_{\{M\}}\left(\Big(-\frac{1}{k},\frac{1}{k}\Big)
\right),
\end{equation*}
and
\begin{equation*}\label{beurlinggerm}
\mathcal{E}_{(M)}^{0}:=\underset{k\in\NN_{>0}}{\varinjlim}\mathcal{E}_{(M)}\left(\Big(-\frac{1}{k},\frac{1}{k}\Big)\right).
\end{equation*}
\end{definition}

Again, if one considers the sequence $M=(j!)_{j \in \NN}$ in the Roumieu case, we obtain the space of {\itshape germs of real analytic functions at} $0\in\RR$; it is denoted by $\mathcal{O}^{0}$.\vspace{6pt}

Let $M\in\RR_{>0}^{\NN}$ be arbitrary and define the sets of weighted formal power series by
$$
\mathcal{F}_{\{M\}}:=\left\{\mathbf{F}=\sum_{j=0}^{+ \infty} F_j x^j
  \, : \, (F_j)_j\in\CC^{\NN}\text{ and }\exists\,h>0 \text{ such that }|\mathbf{F}|^M_h<+\infty\right\},
$$
$$\mathcal{F}_{(M)}:=\left\{\mathbf{F}=\sum_{j=0}^{+ \infty} F_j x^j\,
  : \, (F_j)_j\in\CC^{\NN}\text{ and }\forall\,h>0\, ,\,|\mathbf{F}|^M_h< + \infty\right\},
$$
with
$$|\mathbf{F}|^M_h:=\sup_{j \in \NN} \frac{|F_j| j !}{h^{j} M_j}=\sup_{j \in \NN} \frac{|F_j|}{h^{j} m_j}.$$
We endow these spaces with their natural topology: $\mathcal{F}_{\{M\}
}$ is an (LB)-space and $\mathcal{F}_{(M) }$ a Fr\'echet space.
Naturally, on $\mathcal{F}_{[M]}$ the addition is defined pointwise by
$$
\mathbf{F} + \mathbf{G} =\sum_{j=0}^{+ \infty} \big(F_{j}+ G_{j} \big) x^j
$$
and the scalar multiplication by
$$
\alpha \mathbf{F}  =\sum_{j=0}^{+ \infty} \alpha F_{j} x^j .
$$

\begin{remark}\label{onetoonepowerseries}
It is clear (e.g. see \cite[Remark 2.1.5]{Borelmapgenericity} for some
explanations) that there does exist a one-to-one correspondence
between $\mathcal{F}_{[M]}$ and $\Lambda^1_{[M]}$, the sequence space
has been introduced in \cite[Def. 2.1.4]{Borelmapgenericity}, by
identifying the coefficients $(F_j)_j$ with a sequence (of complex
numbers). So all results from \cite{Borelmapgenericity} (and from \cite{borelmappingquasianalytic}) are also valid
for the sets $\mathcal{F}_{[M]}$ instead of $\Lambda^1_{[M]}$. Note
that in \cite{Borelmapgenericity} we have preferred to work with
classes $\Lambda^1_{[M]}$, but in this present work it seems to be
more natural to consider instead classes of weighted formal power
series as defined above since the Cauchy product $\ast$ seems to be
more natural when considered on $\mathcal{F}_{[M]}$. Note however that we will also obtain results using the pointwise product.
\end{remark}

We introduce the Borel map $j^{\infty}$ (at $0$) by setting
\begin{equation*}\label{Borelmapsequence}
j^{\infty}:\mathcal{E}^{0}_{[M]}\longrightarrow\mathcal{F}_{[M]},\hspace{20pt}j^\infty(f) = \sum_{j=0}^{+ \infty}\frac{f^{(j)}(0)}{j!}x^j.
\end{equation*}

We consider the following definition, according to \cite[Section 2.2]{borelmappingquasianalytic} and \cite[Definition 2.2.1]{Borelmapgenericity}.

\begin{definition}\label{def_weightseq}
A sequence of positive real numbers $M=(M_j)_{j \in \NN} \in\RR_{>0}^{\NN}$ is called a \emph{weight sequence} if
\begin{itemize}
\item[$(I)$] $1=M_0\le M_1$ (normalization),
\item[$(II)$] $M$ is log-convex,
\item[$(III)$] $\liminf_{j\rightarrow\infty}(m_j)^{1/j}>0$.
\end{itemize}
Recall that $m_j:=\frac{M_j}{j!}$ for every $j \in \NN$.
\end{definition}

If $M$ is log-convex and normalized, then $M$ and $j\mapsto(M_j)^{1/j}$ are both increasing and $M_jM_k\le M_{j+k}$ holds for all $j,k\in\NN$, e.g. see \cite[Lemmata 2.0.4, 2.0.6]{diploma}.

Occasionally, we will also consider sequences belonging to the set
$$\hypertarget{LCset}{\mathcal{LC}}:=\{M\in\RR_{>0}^{\NN}:\,M\,\text{normalized, log-convex},\,\lim_{k\rightarrow+\infty}(M_k)^{1/k}=+\infty\}.$$
So for any $M\in\hyperlink{LCset}{\mathcal{LC}}$, assumption $(III)$ above is not necessarily required.

\medskip

Let us also introduce some  classical conditions on a sequence $M\in\RR_{>0}^{\NN}$:
\begin{itemize}
\item $M$ has {\itshape moderate growth}, denoted by \hypertarget{mg}{$(\text{mg})$}, if
$$\exists\,C\ge 1\,\forall\,j,k\in\NN:\,M_{j+k}\le C^{j+k} M_j M_k.$$
\item $M$ is called {\itshape non-quasianalytic,} denoted by \hypertarget{mnq}{$(\text{nq})$}, if
$$\sum_{j=1}^{\infty}\frac{M_{j-1}}{M_j}<+\infty.$$
If $M$ is log-convex, then using {\itshape Carleman's inequality} one
can show (for a proof see e.g. \cite[Proposition 4.1.7]{diploma}) that
$\sum_{j=1}^{\infty}\frac{M_{j-1}}{M_j}<+\infty\Leftrightarrow\sum_{j=1}^{\infty}\frac{1}{(M_j)^{1/j}}<+\infty$.
\item  $M$  is {\itshape quasianalytic} if it does not satify (nq).
\end{itemize}

Let us recall the following classical definition.
\begin{definition}
A subclass $\mathcal{Q}\subseteq\mathcal{E}$ is called {\itshape quasianalytic} if for any open connected set $U\subseteq\RR$ and each point $a\in U$, the Borel map at $a$, denoted by $j^{\infty}_a$, is {\itshape injective} on $\mathcal{Q}(U)$.
\end{definition}

In the case $\mathcal{Q}\equiv\mathcal{E}_{[M]}$ the {\itshape
  Denjoy-Carleman theorem} characterizes this behavior in terms of the
defining weight sequence $M$. More precisely, it states that
$\mathcal{E}_{[M]}$ is quasianalytic if and only if $M$ does not
satisfy \hyperlink{mnq}{$(\text{nq})$}. Let us moreover mention  that
$\mathcal{E}_{[M]}$ is quasianalytic if and only if there do not exist
non-trivial functions in $\mathcal{E}_{[M]}$ with compact support,
e.g. see \cite[Thm. 19.10]{rudin}.
Functions in quasianalytic classes can be represented via a summation
method, obtained within the first part of the proof of \cite[Theorem
3]{thilliez}.

\begin{theorem}[Representation formula, \cite{thilliez}]\label{prop_representationformula}
Let $M$ be a quasianalytic weight sequence. There exist numbers $(\omega^M_{j,k})_{j,k \in \NN}$ such that
\begin{equation*}\label{eq:lim_omega}
\lim_{k \to + \infty} \omega^M_{j,k} =1,  \quad \forall j \in \NN,
\end{equation*}
and such that, given any function $f \in \mathcal{E}^{0}_{\{M\}}$, one has
\begin{equation*}\label{representation formula}
f(x) = \lim_{k \to + \infty} \sum_{j=0}^{k-1} \omega^M_{j,k} \frac{f^{(j)} (0)}{j !} x^j
\end{equation*}
for every $x>0$ small enough.
\end{theorem}

Keeping the notations of this Theorem, we directly obtained in
\cite[Corollary 3.1.2]{Borelmapgenericity} the following important
result. It will be the key for the proofs of algebrability.

\begin{corollary}\label{cor_Thilliez}
Let $M$ be a quasianalytic weight sequence. If
$\mathbf{F}=\sum_{j=0}^{+ \infty}F_{j}x^{j} $ is a formal power series
 for which there exists a sequence of positive real numbers $(a_n)_{n \in \NN}$ decreasing to $0$ such that
\begin{equation*}\label{eq_cor_Thilliez}
\limsup_{k \to + \infty} \left| \sum_{j=0}^{k-1} \omega^{M}_{j,k} F_j a_n^j \right|= + \infty
\end{equation*}
for all $n \in \NN$, then $\mathbf{F} \notin j^\infty (\mathcal{E}^{0}_{\{M\}})$.
\end{corollary}

Finally, let us recall some relations between weight sequences. For two weight sequences $M=(M_j)_j$ and $N=(N_j)_j$ we write $M\le N$ if and only if $M_j\le N_j\Leftrightarrow m_j\le n_j$ holds for all $j\in\NN$. 
Moreover we define $M\hypertarget{mpreceq}{\preceq}N$ by
$$\,\exists\,h,C>0\text{ such that } \forall\,j\in\NN\, , \,\,  M_j\le C h^j N_j$$
or equivalently
$$\sup_{j\in\NN_{>0}}\left(\frac{M_j}{N_j}\right)^{1/j}<+\infty.$$
We call the weight sequences $M$ and $N$  equivalent, denoted by $M\hypertarget{approx}{\approx}N$, if
$$M\hyperlink{mpreceq}{\preceq}N\,\text{and}\,N\hyperlink{mpreceq}{\preceq}M.$$
Finally, we write $M\hypertarget{mtriangle}{\vartriangleleft}N$
if $$\forall\,h>0 \;\exists\,C>0 \text{ such that } \forall \, j\in\NN
\, , \,\,  M_j\le C h^j N_j$$
which is equivalent to
$$
\lim_{j\rightarrow\infty}\left(\frac{M_j}{N_j}\right)^{1/j}=0.
$$
\vspace{6pt}
In the relations above one can replace $M$ and $N$ simultaneously by $m$ and $n$ because the factorial term is cancelling out.\vspace{6pt}

Those relations between weight sequences imply inclusions between ultradifferentiable classes, see  e.g.  \cite[Section 2.2]{borelmappingquasianalytic} and the references therein. More precisely, let $M$ be a weight sequence and $N$ arbitrary, then $M\hyperlink{mpreceq}{\preceq}N$ if and only if $\mathcal{E}_{[M]}\subseteq\mathcal{E}_{[N]}$, which is equivalent to $\mathcal{F}_{[M]}\subseteq\mathcal{F}_{[N]}$. In particular, choosing $M=(j!)_{j\in\NN}$, we get $\mathcal{C}^{\omega}\subseteq\mathcal{E}_{\{N\}}$ if and only if $\liminf_{j\rightarrow + \infty}(n_j)^{1/j}>0$. Moreover, if $N$ is a weight sequence, then $\mathcal{E}_{\{N\}}\subseteq\mathcal{C}^{\omega}$ if and only if $\sup_{j\in\NN_{>0}}(n_j)^{1/j}<+\infty$. Hence $\mathcal{C}^{\omega}\subsetneq\mathcal{E}_{\{N\}}$ if and only if $\sup_{j\in\NN_{>0}}(n_j)^{1/j}=+\infty$.\vspace{6pt}

Similarly $M\hyperlink{mtriangle}{\vartriangleleft}N$  if and only if $\mathcal{E}_{\{M\}}\subsetneq\mathcal{E}_{(N)}$, which is equivalent to $\mathcal{F}_{\{M\}}\subsetneq\mathcal{F}_{(N)}$. In particular, $\mathcal{C}^{\omega}\subsetneq\mathcal{E}_{(N)}$ if and only if $\lim_{j\rightarrow + \infty}(n_j)^{1/j}=+\infty$.

\medskip
Let us close this section by gathering some comments from \cite{Borelmapgenericity}.

\begin{itemize}
\item
In the following sections we will study the Borel map $j^{\infty}$ defined on quasianalytic ultradifferentiable classes such that $\mathcal{C}^{\omega}\subsetneq\mathcal{E}_{[M]}$ holds true. The general assumptions $(I)-(III)$ on $M$ are not restricting the generality of our considerations: For any $M\in\RR_{>0}^{\NN}$ with $\mathcal{C}^{\omega}\subseteq\mathcal{E}_{[M]}$ we have $\liminf_{j\rightarrow + \infty}(m_j)^{1/j}>0$ in the Roumieu and $\lim_{j\rightarrow+ \infty}(m_j)^{1/j}=+\infty$ in the Beurling case and we can replace $M$ by its log-convex minorant $M^{\on{lc}}$ (see \cite[Chapitre I]{mandelbrojtbook} and \cite[(3.2)]{Komatsu73}) without changing the associated ultradifferentiable class whereas only $\mathcal{F}_{[M^{\on{lc}}]}\subseteq\mathcal{F}_{[M]}$ follows (and the weight matrix/function setting is reduced to the sequence case situation).
\item In this paper all the spaces and results are considered in $\RR$,
  but everything goes similarly in $\RR^{r}$ by using a simple
  reduction argument.
\item Finally by translation all results below also hold true if $0\in\RR$ is replaced
by any other point $a\in\RR$.
\end{itemize}

\section{Algebrability with respect to the Cauchy product}\label{Algebrconvolution}

The classical product that can be considered on the space  $\mathcal{F}_{[M]}$  is the Cauchy product (or convolution). It is defined by
$$\mathbf{F} \ast \mathbf{G} =\sum_{j=0}^{+ \infty} \left( \sum_{r=0}^{j} F_{r} G_{j-r}\right) x^j.$$
The aim of this section is to obtain results of algebrability in
$\mathcal{F}_{[M]}$ endowed with the Cauchy product. Then we extend them to the weight matrix and weight function settings.

By the Leibnitz formula, we have that pointwise multiplication of
functions is transferred to the Cauchy product for their formal power
series, i.e. one has $j^{\infty}(fg) = j^{\infty}(f)\ast
j^{\infty}(g)$. A proof for the closedness under the pointwise product
of ultradifferentiable functions can be found in \cite[Proposition
2.0.8]{diploma}. By repeating these arguments we can show the following result which ensures that under relatively weak assumptions on $M$ it makes sense to consider the question of algebrability on $\mathcal{F}_{[M]}$.

\begin{lemma}\label{closednessunderconvolution}
If $M=(M_j)_j$ satisfies
\begin{equation}\label{algebracondition}
\exists\,C\ge 1\text{ such that }\forall\,j,k\in\NN \, , \,\,M_jM_k\le C^{j+k}M_{j+k},
\end{equation}
which is the case if $M$ is a normalized
log-convex sequence (see \cite[Lemma 2.0.6]{diploma}), then $\mathcal{F}_{[M]}$ is a ring under $\ast$.
\end{lemma}

\demo{Proof}
Indeed, if
$$
|F_j|\leq \frac{C_{1} h_{1}^{j} M_j}{j!} \,  , \,\, \forall j \in \NN \quad \text{ and } \quad |G_j|\leq
\frac{C_{2} h_{2}^{j} M_j }{j!} \,  , \,\, \forall j \in \NN
$$
for some $C_{1},C_{2},h_{1},h_{2} >0$, then one has
$$
 \left| \sum_{r=0}^{j} F_{r} G_{j-r}\right|
 \le  C_{1}C_{2} \sum_{r=0}^{j}\frac{ h_{1}^{r} M_r}{r!}
        \frac{h_{2}^{j-r} M_{j-r}}{(j-r)!}
 \le  C_{1}C_{2} C^{j}M_{j} \sum_{r=0}^{j}\frac{ h_{1}^{r}}{r!}
        \frac{h_{2}^{j-r} }{(j-r)!} = \frac{C_{3} h^{j}
          M_{j}}{j!}
$$where $C_{3}= C_{1}C_{2}$ and $h= C(h_{1}+h_{2})$.

\qed\enddemo



\subsection{The weight sequence setting}\label{infinityalgebraweightsequence}
We start with the single weight sequence case and  prove the following result.

\begin{theorem}\label{thm:algebraabinfinity}
Let $M$ and $N$ be two quasianalytic weight sequences. Assume that $\mathcal{O}^{0}\subsetneq\mathcal{E}^{0}_{(N)}$ resp. $\mathcal{O}^{0}\subsetneq\mathcal{E}^{0}_{\{N\}}$, i.e.
\begin{equation}\label{strictincl}
\lim_{j\rightarrow +\infty}(n_j)^\frac{1}{j}=+\infty\hspace{30pt}\text{resp.}\hspace{30pt}\sup_{j\in\NN_{>0}}\left(n_{j}\right)^\frac{1}{j}=+\infty.
\end{equation}

Then $\mathcal{F}_{[N]} \setminus j^{\infty}(\mathcal{E}^{0}_{\{M\}})$ is $\mathfrak{c}$-algebrable in $\mathcal{F}_{[N]}$ endowed with the Cauchy product (hence $\mathcal{F}_{[N]} \setminus j^{\infty}(\mathcal{E}^{0}_{(M)})$ too).
\end{theorem}

\demo{Proof}
By assumption, we can consider an increasing sequence $(k_p)_{p \in \NN} $ of natural numbers satisfying:
\begin{enumerate}[(i)]
\item $k_0=1$ and $k_p > p k_{p-1}$ for every $p \in \NN_{>0}$,
\item 
$\lim_{p\rightarrow + \infty}\left(n_{k_p}\right)^\frac{1}{k_p}=+\infty$, 
\item $\sum_{j=0}^{pk_{p-1}}\left|\omega^M_{j,k_p}-1\right|n_j\le 1$
  for every $p \in \NN_{>0}$, where the numbers
  $(\omega^M_{j,k})_{j,k\in\NN}$ are those arising in Theorem \ref{prop_representationformula}.
\end{enumerate}

Let $(A,B)$ be an open interval with $0<A<B<1$. Let us also consider  a {\itshape Hamel basis} $\mathcal{H}$ of $\RR$ (i.e. a basis of $\RR$ seen as a $\QQ$ vector space). We can assume that the elements of $\mathcal{H}$ are in $(A,B)$. Indeed, if $h \in \mathcal{H}$ is not in $(A,B)$, it suffices to consider $q_h  \in \QQ$ such that $q_h h \in (A,B)$, and we keep a basis.

 For an arbitrary given value $b\in\mathcal{H}$, we define the formal power series $\mathbf{F}^b = \sum_{j=0}^{+ \infty} F^b_j x^j$  by setting
\begin{equation*}\label{squenceF}
F^b_j:=\left\{
 \begin{array}{cl}
 (n_{k_{p}})^{b}&\,\text{ if }\,j= k_{p}, \\[1.5ex]
 0 & \,\text{ if } j \notin \{k_{p}: p \in \NN \}.
 \end{array}
\right.
\end{equation*}
Since $b<1$, it is straightforward to check that $\mathbf{F}^b\in
\mathcal{F}_{(N)}$ (hence also $\mathbf{F}^b\in\mathcal{F}_{\{N\}}$). 
Let us note that if $\mathbf{F}:=\mathbf{F}^b$ and if we define  the formal power series $\mathbf{F}^{(i)}:= \underbrace{\mathbf{F} \ast \dots \ast \mathbf{F}}_{i \text{ times}}$, then one has
\begin{equation*}\label{eq:puissF}
F^{(i)}_j = \sum_{(k_{p_1},\dots,k_{p_i}) \in \mathcal{A}^i(j)}
F_{k_{p_1}} \cdots F_{k_{p_i}} = \sum_{(k_{p_1},\dots,k_{p_i}) \in
  \mathcal{A}^i(j)} \big(n_{k_{p_1}} \cdots n_{k_{p_i}} \big)^{b},
\end{equation*}
where
$$\mathcal{A}^i(j) := \big\{(k_{p_1},\dots,k_{p_i}) \in \mathcal{A}^i : k_{p_1} + \dots + k_{p_i} = j \big\}.
$$
In particular, if $j = ik_p$ for some $p \geq i$ and if $k_{p_1}+
\dots+k_{p_i} = ik_p$, then one has  $k_{p_1}= \dots =k_{p_i} = k_p$
since the sequence $(k_q)_{q \in \NN}$ is strictly increasing and
since $k_{p+1}> (p+1) k_{p} > i k_{p}$. Consequently, one has
$\mathcal{A}^i(ik_{p})  = \{ (k_{p}, \dots, k_{p})\}$ and
$$
F^{i}_{ik_{p}} = \big(n_{k_{p}}\big)^{ib}
$$
if $p \geq i$. Note also that $F^{(i)}_j = 0$ if $j \in \{pk_{p-1}+1,
\dots, k_{p}-1\}$ if $p \geq i$ since in this case $\mathcal{A}^i(j) = \emptyset$.

\medskip
Now, let us consider the algebra $\mathcal{G}$ generated by $\big\{\mathbf{F}^b: b\in\mathcal{H}\big\}$ and let us show that  $\mathcal{G}$  has the desired property. Any element of this algebra can be written as
$$
\mathbf{G} = \sum_{l=1}^L \alpha_l \underbrace{\left(\mathbf{F}^{b_1} \ast \cdots \ast \mathbf{F}^{b_1}\right)}_{i_{l,1} \text{ times}}\ast\cdots\ast\underbrace{\left(\mathbf{F}^{b_J} \ast \cdots \ast \mathbf{F}^{b_J}\right)}_{i_{l,J} \text{ times}},
$$
where $\alpha_1, \dots, \alpha_L \neq 0$, $b_1, \dots, b_{J}\in \mathcal{H}$ are pairwise distinct, for every $l\in \{1, \dots, L\}$ there is at least one $m \in \{1,\dots, J\} $ such that $i_{l,m} \neq 0$ and for every $l, l' \in \{1, \dots, L\}$, $l \neq l'$, there is at least one $m \in \{1,\dots, J\}$ such that $i_{l,m} \neq i_{l',m}$. For every $l \in \{1, \dots, L\}$, let us set
$$
P_l := i_{l,1} + \dots + i_{l,J}.
$$

As done in the case of a single power series, if $p \geq P_l$, one has
\begin{equation}\label{puissFgeneral}
\Bigg(\underbrace{\left(\mathbf{F}^{b_1} \ast \cdots \ast \mathbf{F}^{b_1}\right)}_{i_{l,1} \text{ times}}\ast\cdots\ast\underbrace{\left(\mathbf{F}^{b_J} \ast \cdots \ast \mathbf{F}^{b_J}\right)}_{i_{l,J} \text{ times}} \Bigg)_{P_l k_p} = (n_{k_p})^{i_{l,1}{b_1}+ \dots + {i_{l,J}}{b_J}}
\end{equation}
and if furthermore $j \in \{pk_{p-1}+1,\dots, k_p - 1\}$, then
$$
\Bigg(\underbrace{\left(\mathbf{F}^{b_1} \ast \cdots \ast \mathbf{F}^{b_1}\right)}_{i_{l,1} \text{ times}}\ast\cdots\ast\underbrace{\left(\mathbf{F}^{b_J} \ast \cdots \ast \mathbf{F}^{b_J}\right)}_{i_{l,J} \text{ times}} \Bigg)_{j} = 0.
$$
It follows that
\begin{equation}\label{puissFgeneral2}
G_j = 0, \quad \forall j \in  \{pk_{p-1}+1,\dots, k_p - 1\},
\end{equation}
as soon as $p \geq P :=\max_{l \in \{1, \dots,L\}}P_l$.

\medskip
In order to show that the formal power series $\mathbf{G}$ does not
belong to the image $j^{\infty}(\mathcal{E}^{0}_{\{M\}})$ of the Borel map, by Corollary \ref{cor_Thilliez} it suffices to show that
$$
\limsup_{k \to + \infty} \left| \sum_{j=0}^{k-1} \omega^{M}_{j,k}G_j a^j \right|= + \infty
$$
for every $a>0$ small enough. Of course, it suffices to prove that
\begin{equation*}
\limsup_{p \to + \infty} \left| \sum_{j=0}^{k_p-1} \omega^{M}_{j,k_p}G_j a^j \right|= + \infty.
\end{equation*}
If $p \geq P$, then by \eqref{puissFgeneral2}, one has
\begin{equation}\label{Gsumgeneral}
\sum_{j=0}^{k_p-1} \omega^{M}_{j,k_p}G_j a^j=\sum_{j=0}^{pk_{p-1}} \omega^{M}_{j,k_p}G_j a^j=\sum_{j=0}^{pk_{p-1}} G_j a^j  + \sum_{j=0}^{pk_{p-1}} \big(\omega^{M}_{j,k_p} - 1\big)G_j a^j .
\end{equation}
The first term of the sum is a power series, so its convergence or divergence properties are easy to study. So, let us start with this expression. We have
\begin{align*}
\limsup_{j \to + \infty} |G_j|^{1/j}& \geq \limsup_{p \to + \infty} \big|G_{P k_p}\big|^{1/(Pk_p)} = \limsup_{p \to + \infty} \Big|\sum_{l : P_l = P} \alpha_l (n_{k_p})^{{i_{l,1}}{b_1}+ \dots + {i_{l,J}}{b_J}}\Big|^{1/(Pk_p)} .
\end{align*}
Note that the exponents ${i_{l,1}}{b_1}+ \dots + {i_{l,J}}{b_J}$ are pairwise distinct. Indeed since the $i_{l,j}$ are natural numbers and since if $l \neq l'$ there is at least one number $j$ such that $i_{l,j} \neq i_{l',j}$, it is impossible to have
$${i_{l,1}}{b_1}+ \dots + {i_{l,J}}{b_J} = {i_{l',1}}{b_1}+ \dots + {i_{l',J}}{b_J},$$
because this would contradict the linear independence of the values ${b_1},\dots, {b_J} \in \mathcal{H}$.
Hence, the desired behavior will be given by the largest one (since $n_{k_p}\rightarrow+\infty$ as $p\rightarrow+\infty$) and we can write
\begin{align*}
\limsup_{p \to + \infty} \Big|\sum_{l : P_l = P} \alpha_l (n_{k_p})^{{i_{l,1}}{b_1}+ \dots + {i_{l,J}}{b_J}}\Big|^{1/(Pk_p)} \geq C \times\limsup_{p \to + \infty} \Big|\alpha_l (n_{k_p})^{{i_{l,1}}{b_1}+ \dots + {i_{l,J}}{b_J}}\Big|^{1/(Pk_p)},
\end{align*}
for some positive constant $C$ and some $l$ well chosen such that $P=P_{l}$. This last expression can be estimated by
\begin{align*}
&\limsup_{p \to + \infty} \Big|\alpha_l\Big|^{1/(Pk_p)}\Big|(n_{k_p})^{{i_{l,1}}{b_1}+ \dots + {i_{l,J}}{b_J}}\Big|^{1/(Pk_p)}\ge\limsup_{p \to + \infty} \Big|\alpha_l\Big|^{1/(Pk_p)}\Big|(n_{k_p})^{A({i_{l,1}}+ \dots + i_{l,J})}\Big|^{1/(Pk_p)}
\\&
=\limsup_{p \to + \infty} \Big|\alpha_l\Big|^{1/(Pk_p)}\Big|n_{k_p}\Big|^{A/k_p}=+\infty,
\end{align*}
by recalling ${i_{l,1}+\dots+i_{l,J}}=P$ and assumption $(ii)$ from above. 
Hence this first term of the sum in \eqref{Gsumgeneral} cannot be bounded.

\medskip
Let us now study the second term of the sum in
\eqref{Gsumgeneral}. Since $\mathcal{F}_{\{N\}}$ is an algebra for the
Cauchy product, we know that $\mathbf{G}\in \mathcal{F}_{\{N\}}$.  
So there exist $h,C>0$ such that
$$
\sup_{j \in \NN} \frac{|G_j|}{h^{j} n_j} < C.
$$
Using  assumption $(iii)$, we obtain
\begin{eqnarray*}
\left|\sum_{j=0}^{p k_{p-1}} \big(\omega^{M}_{j,k_p} - 1\big)G_j a^j \right|
& = & \sum_{j=0}^{pk_{p-1}} \left|\omega^{M}_{j,k_p} - 1 \right| |G_j| a^j  \\
& \leq & C\sum_{j=0}^{pk_{p-1}} \left|\omega^{M}_{j,k_p} - 1 \right| n_j (ha)^j  \\
& \leq & C \sum_{j=0}^{pk_{p-1}} \left|\omega^{M}_{j,k_p} - 1 \right| n_j   \\
& \leq & C
\end{eqnarray*}
if $p \geq P$ and $a< \frac{1}{h}$. The conclusion follows.


\qed\enddemo

We wish to mention that each algebra contained in $\mathcal{F}_{[N]} \setminus j^{\infty}(\mathcal{E}^{0}_{\{M\}})$, hence in particular the algebra $\mathcal{G}$ constructed in the previous result, does not contain the identity $\mathbf{E}=1$ for the convolution $\ast$ anymore. Here $E_j=\delta_{j,0}$ and clearly $\mathbf{E}=j^{\infty}(1)$ with $1: x\mapsto 1$ for all $x\in\RR$. Similarly this will be the case for the weight matrix and weight function case below as well.

\subsection{The general weight matrix case}\label{generalmatrixcase}

The aim of this subsection is to establish an equivalent of Theorem
\ref{thm:algebraabinfinity} in the more general setting supplied by weight
matrices.
First we recall the definitions given in \cite[Section 4.1]{Borelmapgenericity}, see also the literature citations therein.

\begin{definition}\label{defsect41}
A {\emph weight matrix} $\mathcal{M}$ is a family of sequences $\mathcal{M}:=\big\{M^{(\lambda)}\in\RR_{>0}^{\NN}: \lambda>0\big\}$, such that
$$\forall\,\lambda >0,\,\,M^{(\lambda)} \text{ is a weight sequence }$$
and
$$\,M^{(\lambda)}\le M^{(\kappa)} 
\text{ for all } \, 0<\lambda\le\kappa.
$$

A matrix is called {\emph constant} if
$M^{(\lambda)}\hyperlink{approx}{\approx}M^{(\kappa)}$
for all $\lambda ,\kappa >0$.
\end{definition}



We introduce classes of ultradifferentiable function of Roumieu type $\mathcal{E}_{\{\mathcal{M}\}}$ and of Beurling type $\mathcal{E}_{(\mathcal{M})}$ as follows (only the pointwise order in Def. \ref{defsect41} is required), see \cite[Section 7]{dissertation} and \cite[Section 4.2]{compositionpaper}.

\begin{definition}
Let $\mathcal{M}$ be a weight matrix and let $U\subseteq\RR$ be non-empty and open. The \emph{$\mathcal{M}$-ultradifferentiable classes of Roumieu and Beurling types} are defined respectively by
\begin{equation*}\label{generalroumieu}
\mathcal{E}_{\{\mathcal{M}\}}(U):=\bigcap_{K\subseteq
  U}\bigcup_{\lambda >0}\mathcal{E}_{\{M^{(\lambda)}\}}(K)
\end{equation*}
and
\begin{equation*}\label{generalbeurling}
\mathcal{E}_{(\mathcal{M})}(U):=\bigcap_{\lambda >0}\mathcal{E}_{(M^{(\lambda)})}(U).
\end{equation*}
\end{definition}

For a compact set $K\subseteq\RR$, one has the representations
$$\mathcal{E}_{\{\mathcal{M}\}}(K):=\underset{\lambda >0}{\varinjlim}\,\underset{h>0}{\varinjlim}\,\mathcal{E}_{M^{(\lambda)},h}(K)$$
and so for $U\subseteq\RR$ non-empty open
\begin{equation*}\label{generalroumieu1}
\mathcal{E}_{\{\mathcal{M}\}}(U)=\underset{K\subseteq
  U}{\varprojlim}\,\underset{\lambda >0}{\varinjlim}\,\underset{h>0}{\varinjlim}\,\mathcal{E}_{M^{(\lambda)},h}(K).
\end{equation*}
Similarly we get for the Beurling case
\begin{equation*}\label{generalbeurling1}
\mathcal{E}_{(\mathcal{M})}(U)=\underset{K\subseteq
  U}{\varprojlim}\,\underset{\lambda >0}{\varprojlim}\,\underset{h>0}{\varprojlim}\,\mathcal{E}_{M^{(\lambda)},h}(K).
\end{equation*}
Consequently, since the sequences of $\mathcal{M}$ are pointwise ordered, $\mathcal{E}_{(\mathcal{M})}(U)$ is a {\itshape Fr\'echet space} and

$\underset{\lambda >0}{\varinjlim}\,\underset{h>0}{\varinjlim}\,\mathcal{E}_{M^{(\lambda)},h}(K)=\underset{n\in\NN_{>0}}{\varinjlim}\,\,\mathcal{E}_{M^{(n)},n}(K)$ is a {\itshape Silva space}, i.e. a countable inductive limit of Banach spaces with compact connecting mappings. For more details concerning the locally convex topology in this setting we refer to \cite[Section 4.2]{compositionpaper}.\vspace{6pt}

\begin{definition}
The spaces of \emph{germs at $0\in\RR$ of the $(\mathcal{M})$-ultradifferentiable functions of Roumieu and Beurling types} are defined respectively by
\begin{equation*}\label{roumieumatrixgerm}
\mathcal{E}_{\{\mathcal{M}\}}^{0}:=\underset{k\in\NN_{>0}}{\varinjlim}\mathcal{E}_{\{\mathcal{M}\}}\left(\Big(-\frac{1}{k},\frac{1}{k}\Big)\right),
\end{equation*}
and
\begin{equation*}\label{beurlingmatrixgerm}
\mathcal{E}_{(\mathcal{M})}^{0}:=\underset{k\in\NN_{>0}}{\varinjlim}\mathcal{E}_{(\mathcal{M})}\left(\Big(-\frac{1}{k},\frac{1}{k}\Big)\right).
\end{equation*}
\end{definition}

Finally, as done in the case of weight sequences, we introduce the corresponding spaces of weighted power series sequences, and we endow them with their classical topology:

$$
\mathcal{F}_{\{\mathcal{M}\}}:=\left\{\mathbf{F}=\sum_{j=0}^{+ \infty}
  F_j x^j \, :
  \,(F_j)_j\in\CC^{\NN}\text{ and }
  \exists  \, \lambda >0 \,, \,\exists \,h>0 \text{ such that }|\mathbf{F}|^{M^{(\lambda)}}_h< + \infty\right\},
$$
$$\mathcal{F}_{(\mathcal{M})}:=\left\{\mathbf{F}=\sum_{j=0}^{+ \infty}
  F_j x^j \, :
  \,(F_j)_j\in\CC^{\NN}\text{ and }
  \forall  \, \lambda >0\,, \,\forall \,h>0 \, , \,|\mathbf{F}|^{M^{(\lambda)}}_h< + \infty\right\}.
$$

Using notations similar as before, the \emph{Borel map} $j^\infty$ is defined in the weight matrix case by
\begin{equation*}\label{Borelmapmatrix}
j^{\infty}:\mathcal{E}^{0}_{[\mathcal{M}]}\longrightarrow\mathcal{F}_{[\mathcal{M}]},\hspace{20pt}j^\infty(f) = \sum_{j=0}^{+ \infty}\frac{f^{(j)}(0)}{j!}x^j.
\end{equation*}

If $\mathcal{M}=\big\{M^{(\lambda)}: \lambda>0\big\}$ is a weight
matrix, then each $M^{(\lambda)}\in\mathcal{M}$ is log-convex and
normalized, i.e. $(I)$ and $(II)$ in Definition \ref{def_weightseq}
are valid. Consequently each $M^{(\lambda)}$ does satisfy \eqref{algebracondition} and thus the proof of Lemma \ref{closednessunderconvolution} together with the fact that the sequences of $\mathcal{M}$ are pointwise ordered immediately imply that both $\mathcal{F}_{\{\mathcal{M}\}}$ and
$\mathcal{F}_{(\mathcal{M})}$ are rings with respect to the
convolution product $\ast$.\vspace{6pt}

Given two matrices $\mathcal{M}$ and $\mathcal{N}$ we write $\mathcal{M}\hypertarget{Mroumprecsim}{\{\preceq\}}\mathcal{N}$ if
$$\forall\;\lambda>0\;\exists\;\kappa>0:\;\;\;M^{(\lambda)}\hyperlink{mpreceq}{\preceq}N^{(\kappa)},$$
and call them {\itshape Roumieu equivalent}, denoted by $\mathcal{M}\hypertarget{Mroumapprox}{\{\approx\}}\mathcal{N}$, if
$\mathcal{M}\hyperlink{Mroumprecsim}{\{\preceq\}}\mathcal{N}$ and $\mathcal{N}\hyperlink{Mroumprecsim}{\{\preceq\}}\mathcal{M}$.

Analogously we write $\mathcal{M}\hypertarget{Mbeurprecsim}{(\preceq)}\mathcal{N}$ if
$$\forall\;\lambda>0\;\exists\;\kappa>0:\;\;\;M^{(\kappa)}\hyperlink{mpreceq}{\preceq}N^{(\lambda)},$$
and call them {\itshape Beurling equivalent}, denoted by $\mathcal{M}\hypertarget{Mbeurapprox}{(\approx)}\mathcal{N}$, if
$\mathcal{M}\hyperlink{Mbeurprecsim}{(\preceq)}\mathcal{N}$ and $\mathcal{N}\hyperlink{Mbeurprecsim}{(\preceq)}\mathcal{M}$.

We have $\mathcal{M}[\preceq]\mathcal{N}$ if and only if $\mathcal{E}_{[\mathcal{M}]}\subseteq\mathcal{E}_{[\mathcal{N}]}$, see \cite[Prop. 4.6]{compositionpaper}.

\begin{definition}\label{matrixquasianalytic}
A weight matrix $\mathcal{M}$ is called \emph{quasianalytic} if for
all $\lambda >0$ the sequence $M^{(\lambda)}$ is quasianalytic. 
\end{definition}

Given a quasianalytic weight matrix $\mathcal{M}$, both classes $\mathcal{E}_{\{\mathcal{M}\}}$ and $\mathcal{E}_{(\mathcal{M})}$ and all classes $\mathcal{E}_{\{M^{(\lambda)}\}}$ resp. $\mathcal{E}_{(M^{(\lambda)})}$ are quasianalytic, too.

If $\mathcal{M}$ is a quasianalytic weight matrix, then to ensure $\mathcal{O}^{0}\subsetneq\mathcal{E}^{0}_{(\mathcal{M})}$ resp. $\mathcal{O}^{0}\subsetneq\mathcal{E}^{0}_{\{\mathcal{M}\}}$ we assume
\begin{equation}\label{strictincl2}
\forall\,\lambda>0\,\, \lim_{j\rightarrow
  +\infty}(m^{(\lambda)}_j)^\frac{1}{j}=+\infty\quad\text{resp.}\quad
\forall\,\lambda>0 \,\,\sup_{j\in\NN_{>0}}\left(m^{(\lambda)}_{j}\right)^\frac{1}{j}=+\infty.
\end{equation}

Let us now prove the generalization of Theorem
\ref{thm:algebraabinfinity} for the matrix setting. The idea of the
proof is based on the following lemma, which allows to reduce the
general case of two weight matrices $\mathcal{N}$ and $\mathcal{M}$ to
the case of a  weight matrix $\mathcal{N}$ and a single weight
sequence $M$ (analogously as done in \cite[Section
4.2]{Borelmapgenericity}).

\begin{lemma}\label{roumieuauxiliary}
Let $\mathcal{M}=\big\{M^{(\lambda)}: \lambda >0\big\}$ be a quasianalytic weight matrix. Then there exists a quasianalytic weight sequence $L$ satisfying $M^{(\lambda)}\hyperlink{mtriangle}{\vartriangleleft}L$ for all $\lambda>0$, i.e. $\mathcal{E}_{\{\mathcal{M}\}} \subseteq \mathcal{E}_{(L)}$ holds true.
\end{lemma}

The general result can be stated as follows.

\begin{theorem}\label{thm:algebraabinfinityfullmatrix}
Let $\mathcal{M}$ and $\mathcal{N}$ be two quasianalytic weight
matrices. Assume that  $\mathcal{O}^{0}\subsetneq\mathcal{E}^{0}_{(\mathcal{N})}$ resp. $\mathcal{O}^{0}\subsetneq\mathcal{E}^{0}_{\{\mathcal{N}\}}$, i.e.
\begin{equation*}\label{strictinclmatrix}
\forall\,\lambda >0\, , \, \lim_{j\rightarrow
  +\infty}\left(n^{(\lambda)}_j\right)^\frac{1}{j}=+\infty\quad\text{
  resp. }\quad \forall\,\lambda>0 \, , \,\sup_{j\in\NN_{>0}}\left(n^{(\lambda)}_{j}\right)^\frac{1}{j}=+\infty.
\end{equation*}
Then $\mathcal{F}_{[\mathcal{N}]} \setminus
j^{\infty}(\mathcal{E}^{0}_{\{\mathcal{M}\}})$ is
$\mathfrak{c}$-algebrable in $\mathcal{F}_{[\mathcal{N}]}$ endowed
with the Cauchy product (hence $\mathcal{F}_{[\mathcal{N}]} \setminus j^{\infty}(\mathcal{E}^{0}_{(\mathcal{M})})$ too).
\end{theorem}



\demo{Proof}
Using Lemma \ref{roumieuauxiliary}, we can consider a quasianalytic
weight sequence $L$ such that $\mathcal{E}_{\{\mathcal{M}\}} \subseteq
\mathcal{E}_{(L)}$. It suffices now to show that  $\mathcal{F}_{[\mathcal{N}]} \setminus
j^{\infty}(\mathcal{E}^{0}_{(L)})$ is
$\mathfrak{c}$-algebrable. The Roumieu case is a consequence of
Theorem \ref{thm:algebraabinfinity}: indeed, it suffices to fix a
weight sequence $N^{(\lambda_{0})}\in \mathcal{N}$ and use the obvious
inclusion $\mathcal{F}_{\{N^{(\lambda_{0})}\}}\subseteq \mathcal{F}_{\{\mathcal{N}\}}$.
For the Beurling case, we will follow the proof of Theorem
\ref{thm:algebraabinfinity}. First, by induction we can construct an increasing sequence $(k_p)_{p \in \NN} $ of natural numbers satisfying:
\begin{enumerate}[(i)]
\item $k_0=1$ and $k_p > p k_{p-1}$ for every $p \in \NN_{>0}$,
\item $\lim_{p\rightarrow + \infty}\left(n^{(1/(p+1))}_{k_p}\right)^\frac{1}{k_p}=+\infty$,
\item $\sum_{j=0}^{pk_{p-1}}\left|\omega^L_{j,k_p}-1\right|n_j^{(p)}\le 1$ for every $p \in \NN_{>0}$.
\end{enumerate}

Then let us consider an open interval $(A,B)$  with
$0<A<B<1$ and a Hamel basis  $\mathcal{H} \subseteq
(A,B)$ of $\RR$. For an arbitrary given value
$b\in\mathcal{H}$, we define the formal power series $\mathbf{F}^b =
\sum_{j=0}^{+ \infty} F^b_j x^j$  by setting
\begin{equation*}\label{squenceF}
F^b_j:=\left\{
 \begin{array}{cl}
 \left(n^{(1/(p+1))}_{k_{p}}\right)^{b}&\text{ if }j = k_{p}, \\[1.5ex]
 0 & \text{ if } j \notin \{k_p:p\in \NN\}.
 \end{array}
\right.
\end{equation*}
It is straightforward to check that $\mathbf{F}^b\in
\mathcal{F}_{(\mathcal{N})}$ for any $b\in \mathcal{H}$. We follow then the lines of the proof of Theorem \ref{thm:algebraabinfinity} where \eqref{puissFgeneral} turns into
\begin{equation*}
\Bigg(\underbrace{\left(\mathbf{F}^{b_1} \ast \cdots \ast \mathbf{F}^{b_1}\right)}_{i_{l,1} \text{ times}}\ast\cdots\ast\underbrace{\left(\mathbf{F}^{b_J} \ast \cdots \ast \mathbf{F}^{b_J}\right)}_{i_{l,J} \text{ times}} \Bigg)_{P_l k_p} = \left(n_{k_p}^{(1/(p+1))}\right)^{{i_{l,1}}{b_1}+ \dots + {i_{l,J}}{b_J}}
\end{equation*}
as soon as $p \geq P_l$. 
We consider again the splitting \eqref{Gsumgeneral} and proceed for
the first term as in Theorem \ref{thm:algebraabinfinity}. Concerning the estimation of the second term of the sum in \eqref{Gsumgeneral}, since $\mathbf{G}\in \mathcal{F}_{\{\mathcal{N}\}}$ 
there exist an index $\lambda_0>0$ and $h,C>0$ such that $$\sup_{j \in \NN} \frac{|G_j|}{h^{j} n^{(\lambda_0)}_j} < C.$$

It follows that
\begin{eqnarray*}
\left|\sum_{j=0}^{p k_{p-1}} \big(\omega^{L}_{j,k_p} - 1\big)G_j a^j \right|
& = & \sum_{j=0}^{pk_{p-1}} \left|\omega^{L}_{j,k_p} - 1 \right| |G_j| a^j  \\
& \leq & C\sum_{j=0}^{pk_{p-1}} \left|\omega^{L}_{j,k_p} - 1 \right| n^{(\lambda_0)}_j (ha)^j  \\
& \leq & C \sum_{j=0}^{pk_{p-1}} \left|\omega^{L}_{j,k_p} - 1 \right| n^{(p)}_j   \\
& \leq & C
\end{eqnarray*}
if $p \geq\max\{P,\lambda_0\}$ and $a< \frac{1}{h}$, and using
assumption $(iii)$. This concludes the proof.
\qed\enddemo

\subsection{The weight function case}\label{weightfunctioncase}
In this section we will study classes of ultradifferentiable functions
defined using weight functions $\omega$ in the sense of Braun, Meise
and Taylor, see \cite{BraunMeiseTaylor90}. As done in \cite{borelmappingquasianalytic} and \cite{Borelmapgenericity}, we will see that this case can be reduced to the weight matrix situation by using the matrix {\emph associated with} $\omega$. First, let us start by recalling the basic definitions.

\begin{definition}
A function $\omega:[0,+\infty)\rightarrow[0,+\infty)$ is called a \emph{weight function} if
\begin{itemize}
\item[$(i)$] $\omega$ is continuous,
\item[$(ii)$] $\omega$ is increasing,
\item[$(iii)$] $\omega(t)=0$ for all $t\in[0,1]$ (normalization, w.l.o.g.),
\item[$(iv)$] $\lim_{t\rightarrow+\infty}\omega(t)=+\infty$.
\end{itemize}
In this case, we say that $\omega$ has $\hypertarget{om0}{(\omega_0)}$.
\end{definition}

Classical additional conditions can be imposed on the considered weight functions. More precisely, let us define the following conditions:
\begin{itemize}
\item[\hypertarget{om1}{$(\omega_1)$}] $\omega(2t)=O(\omega(t))$ as $t\rightarrow+\infty$,

\item[\hypertarget{om2}{$(\omega_2)$}] $\omega(t)=O(t)$ as $t\rightarrow+ \infty$,

\item[\hypertarget{om3}{$(\omega_3)$}] $\log(t)=o(\omega(t))$ as $t\rightarrow+\infty$ ($\Leftrightarrow\lim_{t\rightarrow+\infty}\frac{t}{\varphi_{\omega}(t)}=0$),

\item[\hypertarget{om4}{$(\omega_4)$}] $\varphi_{\omega}:t\mapsto\omega(e^t)$ is a convex function on $\RR$,

\item[\hypertarget{om5}{$(\omega_5)$}] $\omega(t)=o(t)$ as $t\rightarrow+\infty$.

\end{itemize}

For convenience, we define the set
$$\hypertarget{omset1}{\mathcal{W}}:=\{\omega:[0,+ \infty)\rightarrow[0, + \infty): \omega\,\text{has}\,\hyperlink{om0}{(\omega_0)}, \hyperlink{om1}{(\omega_1)}, \hyperlink{om3}{(\omega_3)},\hyperlink{om4}{(\omega_4)}\}.$$
Note that \hyperlink{om2}{$(\omega_2)$} is sometimes also considered as a general assumption on $\omega$ (e.g. see \cite[Sect. 4.1]{borelmappingquasianalytic}) and note also that \hyperlink{om5}{$(\omega_5)$} implies \hyperlink{om2}{$(\omega_2)$}.

For $\omega\in\hyperlink{omset}{\mathcal{W}}$, we define the {\itshape Legendre-Fenchel-Young-conjugate} of $\varphi_{\omega}$ by
$$\varphi^{*}_{\omega}(x):=\sup\{x y-\varphi_{\omega}(y): y\ge 0\},\,\,\,x\ge 0.$$

\begin{definition}
Let  $U\subseteq\RR$ be a non-empty open set and let $\omega\in\hyperlink{omset}{\mathcal{W}}$. The $\omega$-\emph{ultradifferentiable Roumieu type class} is defined by
$$\mathcal{E}_{\{\omega\}}(U):=\{f\in\mathcal{E}(U):\,\forall\,K\subseteq U\,\text{compact}\,\, \exists\,l>0,\,\|f\|^{\omega}_{K,l}<+\infty\},$$
and the $\omega$-\emph{ultradifferentiable Beurling type class} by
$$\mathcal{E}_{(\omega)}(U):=\{f\in\mathcal{E}(U):\,\forall\,K\subseteq U\,\text{compact}\,\, \forall\,l>0,\,\|f\|^{\omega}_{K,l}<+\infty\},$$
where we have put
\begin{equation*}\label{semi-norm-1}
\|f\|^{\omega}_{K,l}:=\sup_{j\in\NN,x\in K}\frac{|f^{(j)}(x)|}{\exp\big(\frac{1}{l}\varphi^{*}_{\omega}(lj)\big)}.
\end{equation*}
\end{definition}

As done in the previous contexts, these spaces are endowed with their
natural topologies. 
Let $\sigma,\tau$ be weight functions, we write
$\sigma\hypertarget{ompreceq}{\preceq}\tau$ if $\tau(t)=O(\sigma(t))$
as $t\rightarrow+\infty$ and call them equivalent, denoted by
$\sigma\hypertarget{sim}{\sim}\tau$, if
$\sigma\hyperlink{ompreceq}{\preceq}\tau$ and
$\tau\hyperlink{ompreceq}{\preceq}\sigma$. Let
$\tau,\sigma\in\hyperlink{omset}{\mathcal{W}}$, then
$\sigma\hypertarget{sim}{\sim}\tau$ if and only if
$\mathcal{E}_{[\sigma]}=\mathcal{E}_{[\tau]}$, see
\cite[Cor. 5.17]{compositionpaper}.

\medskip

Analogously as in the sections above, we also consider the spaces of
germs at $0$, denoted by $\mathcal{E}_{\{\omega\}}^{0}$ and
$\mathcal{E}_{(\omega)}^{0}$, and the associated spaces of weighted
power series $\mathcal{F}_{\{ \omega\} }$ and $\mathcal{F}_{(\omega)
}$. Again, we endow these spaces  with their natural topology: $\mathcal{F}_{\{ \omega\} }$ is an (LB)-space and $\mathcal{F}_{(\omega) }$ a Fr\'echet space. In this setting, the  Borel map is given by
\begin{equation*}\label{Borelmapfunction}
j^{\infty}:\mathcal{E}_{[\omega]}^{0}\longrightarrow\mathcal{F}_{[\omega]},\hspace{20pt}j^\infty(f)=\mathbf{F}=\sum_{j=0}^{+\infty}\frac{f^{(j)}(0)}{j!}x^j.
\end{equation*}

As pointed out in \cite[Section 4.2]{borelmappingquasianalytic}, that
to ensure $\mathcal{C}^{\omega}\subsetneq\mathcal{E}_{\{\omega\}}$
resp. $\mathcal{C}^{\omega}\subsetneq\mathcal{E}_{(\omega)}$, one has
to assume that
\begin{equation}\label{strictinclusioncomega}
\liminf_{t\rightarrow+
  \infty}\frac{\omega(t)}{t}=0\hspace{30pt}\text{resp.}\hspace{30pt}\omega(t)=o(t)\,\text{
  as }\,t\rightarrow+ \infty ,\,\text{i.e.}\,\hyperlink{om5}{(\omega_5)},
\end{equation}
which follows from the characterizations given in \cite[Lemma 5.16, Cor. 5.17]{compositionpaper} and the fact that the weight $\omega(t)=t$ (up to equivalence) defines the class $\mathcal{C}^{\omega}$.

Moreover, in the present setting, the definition of quasianalyticity takes the following form.

\begin{definition}
A weight function is called \emph{quasianalytic} if it satisfies
\begin{equation}\label{omegaQ}
\int_1^{+ \infty}\frac{\omega(t)}{t^2}dt= + \infty.
\end{equation}
\end{definition}

In \cite{dissertation} and \cite[Section 5]{compositionpaper}, a matrix $\Omega:=\big\{W^{(\lambda)}=(W^{(\lambda)}_j)_{j\in\NN}: \lambda>0\big\}$ has been associated with each $\omega\in\hyperlink{omset}{\mathcal{W}}$: This matrix is defined by
$$W^{(\lambda)}_j:=\exp\left(\frac{1}{\lambda}\varphi^{*}_{\omega}(\lambda j)\right) , \quad \forall j\in\NN, \, \forall\lambda >0,$$
and  $\mathcal{E}_{[\omega]}=\mathcal{E}_{[\Omega]}$ holds as locally convex vector spaces. Moreover, the following results have been obtained (for which \hyperlink{om1}{$(\omega_1)$} is not needed necessarily):

\begin{itemize}
\item[$(i)$] Each $W^{(\lambda)}$ satisfies the basic assumptions $(I)$ and $(II)$ and $\lim_{j\rightarrow+\infty}(W^{(\lambda)}_j)^{1/j}=+\infty$.

\item[$(ii)$] $\omega$ has in addition \hyperlink{om2}{$(\omega_2)$} if and only if some/each $W^{(\lambda)}$ has $(III)$, too.
\end{itemize}

So each $W^{(\lambda)}\in\Omega$ is a weight sequence according to the requirements from Definition \ref{def_weightseq}, provided $\omega\in\hyperlink{omset}{\mathcal{W}}$ has \hyperlink{om2}{$(\omega_2)$}. Moreover, by \cite[Corollary 5.8]{compositionpaper} and \cite[Corollary 4.8]{testfunctioncharacterization}, one has that the following assertions are equivalent (again \hyperlink{om1}{$(\omega_1)$} is not needed but then $\mathcal{E}_{[\omega]}=\mathcal{E}_{[\Omega]}$ fails):
\begin{itemize}
\item[$(i)$] $\omega\in\hyperlink{omset}{\mathcal{W}}$ is quasianalytic,
\item[$(ii)$] $\Omega$ is quasianalytic in the sense of Definition \ref{matrixquasianalytic},
\item[$(iii)$] some/each $W^{(\lambda)}$ is quasianalytic.
\end{itemize}

Similarly, from \cite[Proposition 2]{borelmappingquasianalytic} (and in the same spirit as in \cite[Section 5]{compositionpaper}), for any $\omega\in\hyperlink{omset}{\mathcal{W}}$ one gets $\mathcal{F}_{[\omega]}=\mathcal{F}_{[\Omega]}$ as locally convex spaces, too.

Since each $W^{(\lambda)}$ satisfies \eqref{algebracondition} and the sequences $W^{(\lambda)}$ are pointwise ordered, as already commented in the general weight matrix case above, by following the proof of Lemma \ref{closednessunderconvolution} it is immediate to see that for any $\omega\in\hyperlink{omset}{\mathcal{W}}$ the sets $\mathcal{F}_{\{\omega\}}$ and $\mathcal{F}_{(\omega)}$ are always rings w.r.t. the convolution product.\vspace{6pt}

The weight function approach is again reduced to the more general weight matrix setting by using the weight matrices $\mathcal{N}=\Omega$ and $\mathcal{M}=\Sigma$ associated with $\omega$ and $\sigma$ and Theorem \ref{thm:algebraabinfinityfullmatrix} turns into the following form.

\begin{theorem}\label{thm:algebraabinfinityfullfunction}
Let $\sigma, \omega \in\hyperlink{omset1}{\mathcal{W}}$ be two
quasianalytic weight functions. Assume that $\omega$ satisfies
\hyperlink{om2}{$(\omega_2)$} and $\liminf_{t\rightarrow+
  \infty}\frac{\omega(t)}{t}=0$ in the Roumieu
resp. $\hyperlink{om5}{(\omega_5)}$ in the Beurling case. Then
$\mathcal{F}_{[\omega]} \setminus
j^{\infty}(\mathcal{E}^{0}_{\{\sigma\}})$ is
$\mathfrak{c}$-algebrable in $\mathcal{F}_{[\omega]}$ endowed with the
Cauchy product (hence $\mathcal{F}_{[\omega]} \setminus j^{\infty}(\mathcal{E}^{0}_{(\sigma)})$ too).
\end{theorem}

\section{Algebrability with respect to the pointwise product}\label{Algebrpointwise}

\subsection{Motivation and solid spaces}
Instead of dealing with the Cauchy product $\ast$ on
$\mathcal{F}_{[N]}$, $\mathcal{F}_{[\mathcal{N}]}$ and
$\mathcal{F}_{[\omega]}$, one can also treat the pointwise product, in the literature also known under {\itshape Hadamard product}: Given $\mathbf{F}=\sum_{j=0}^{+\infty}F_jx^j$ and $\mathbf{G}=\sum_{j=0}^{+\infty}G_jx^j$ we consider
\begin{equation}\label{pointwisemulti}
\mathbf{F}\odot\mathbf{G}:=\sum_{j=0}^{+\infty}F_jG_jx^j.
\end{equation}

On the one hand, the study of the problem of algebrability with respect
to this product might be a quite natural question. Moreover this
product has become important very recently by the development of a
convenient theory of multisummability of formal power series, see
\cite[Chapter 4]{dissertationjimenez} and \cite{multisummability}. Concerning these recent insights, in a private communication Prof. J. Sanz has told the authors the following explanations.

\begin{remark}\label{javiscomment}
The natural procedure for assigning a sum to a summable series (in a one step procedure) precisely starts by termwise dividing the coefficients of the series by a moment sequence (equivalent to the weight sequence defining the level) to make the new series (the formal Borel transform) convergent. Correspondingly, the formal Laplace transform multiplies coefficients by the weight sequence. Moreover, sometimes series are not summable but multisummable, i.e. a sum is assigned to them after a finite number of summability procedures, each associated to a different (that is, associated to nonequivalent weight sequences) level, and then one needs to move from one level to another one, which means one has to termwise multiply or divide the coefficients of a given series by a sequence which measures the "jump" between two different levels.
\end{remark}

Consequently, when working within the framework of weight matrices,
one can control these movements/jumps in the sense that one can stay
within a given matrix $\mathcal{M}$ by multiplying pointwise one
sequence $M^1\in\mathcal{M}$ by another one $M^2\in\mathcal{M}$; and
for this behavior closedness under $\odot$ of $\mathcal{F}_{[\mathcal{M}]}$ becomes interesting and crucial.

\medskip

But the study of $\odot$ has also been motivated by the following
approach (cf. \cite{borelmappingquasianalytic},
\cite{Borelmapgenericity}): It is still an open problem to give a
precise characterization which $\mathbf{F}\in\mathcal{F}_{[M]}$ do
belong to the image $j^{\infty}(\mathcal{E}^{0}_{[M]})$ in the
quasianalytic setting (strictly containing the real analytic germs)
and similarly for the weight matrix and weight function setting.
Unlike what happens in the case of Cauchy product, let us show that in this setting
this image $j^{\infty}(\mathcal{E}^{0}_{\{M\}})$ is never closed under
pointwise product $\odot$.

Let us first start by recalling the two following
results. The first one is due to Thilliez, see \cite[Theorem 1]{thilliez} and for a detailed proof also \cite[Prop. 3.1.2]{diploma} and \cite[Lemma 2.9]{compositionpaper}.

\begin{proposition}\label{prop_thetafunction}
Let $M\in \RR^{\NN}_{>0}$ satisfying the conditions $(I)$ and $(II)$
from Definition \ref{def_weightseq}. Let us consider the function
$$\theta_{M}(x)= \sum_{k=0}^{\infty}\frac{M_k}{(2\mu_k)^k}\exp(2i\mu_k
x), \quad x \in \RR,$$
with $\mu_k:=M_k/M_{k-1}$ for $k \in \NN_{>0}$ and $\mu_0:=1$. Then
$\theta_M\in\mathcal{E}_{\{M\}}(\RR)$ and
\begin{equation*}\label{thetafunction}
\theta_M^{(j)}(0)=i^j s_{j}  \text{ with } s_{j}\ge M_{j}\, ,  \quad \forall\,j\in\NN.
\end{equation*}
\end{proposition}

It is not difficult to see that such a function $\theta_M$ does not belong to the Beurling
type class associated to $M$. On the opposite direction, \cite[Thm. 2]{borelmappingquasianalytic} and its proof show that if the derivatives of a smooth function $f$ at $0$ have ''large size'' and all have the same sign, then $f$ cannot
belong to any quasianalytic germ class $\mathcal{E}^{0}_{\{M\}}$. More precisely we have:

\begin{proposition}\label{prop_positive}
Let $M$ be a quasianalytic weight sequence satisfying
$\mathcal{O}^{0}\subsetneq\mathcal{E}^{0}_{\{M\}}$. Assume that the
formal power series
$\mathbf{F}=\sum_{j=0}^{+\infty}F_jx^j\in\mathcal{F}_{\{M\}}$ with
$F_j>0$ for all $j\in\NN$ does not define a real analytic germ. Then
$\mathbf{F} \notin j^{\infty}(\mathcal{E}^{0}_{\{N\}})$ for any quasianalytic weight sequence $N$.
\end{proposition}


These two results lead to the following observation:
If $M$ is a quasianalytic weight sequence such that
$\mathcal{O}^{0}\subsetneq\mathcal{E}^{0}_{\{M\}}$, then there does exist $\mathbf{F}\in\mathcal{F}_{\{M\}}$ such that $\mathbf{F}\in j^{\infty}(\mathcal{E}^{0}_{\{M\}})$ but $|\mathbf{F}|\notin j^{\infty}(\mathcal{E}^{0}_{\{M\}})$ with $|\mathbf{F}|:=\sum_{j=0}^{+\infty}|F_j|x^j$.
Indeed, it suffices to consider $\mathbf{F}= j^{\infty}(\theta_{M})$, i.e. $F_j:=\theta_M^{(j)}(0)/j!$.\vspace{6pt}

{\itshape Conclusion:} Multiplying a given $\mathbf{F}\in j^{\infty}(\mathcal{E}^{0}_{\{M\}})$ pointwise by a formal power series $\mathbf{S}$ given in terms of a sequence of suitable complex numbers on the unit circle, and so $\mathbf{S}\in\mathcal{F}_{\{(j!)_j\}}=j^{\infty}(\mathcal{O}^0)\subseteq j^{\infty}(\mathcal{E}^{0}_{\{M\}})$ is obvious, will in general yield that $\mathbf{F}\odot\mathbf{S}=|\mathbf{F}|\notin j^{\infty}(\mathcal{E}^{0}_{\{M\}})$. Thus for any quasianalytic weight sequence $M$ with $\mathcal{O}^{0}\subsetneq\mathcal{E}^{0}_{\{M\}}$ closedness under $\odot$ fails for the space $j^{\infty}(\mathcal{E}^{0}_{\{M\}})$. Note that $j^{\infty}(\mathcal{E}^{0}_{\{M\}})$ is closed under $\ast$ for any weight sequence $M$ by having $j^{\infty}(fg)=j^{\infty}(f)\ast j^{\infty}(g)$.\vspace{6pt}

Connected to this observation is the notion of {\itshape solid sub-
  and superspaces} for spaces of (complex) sequences, e.g. see
\cite{andersonshields}. Let $A$ be a vector spaces of sequences, then
$A$ is said to be {\itshape solid} if $(a_j)_j\in A$ does imply
$(b_j)_j\in A$ for all sequences satisfying $|b_j|\le|a_j|$,
$\forall j\in\NN$.

\medskip

In \cite[Lemma 2]{andersonshields} it has been shown that for any
given sequence space $A$ there does exist $s(A)$, the {\itshape
  largest solid subspace (or solid core)} of $A$, and there does exist
$S(A)$, the {\itshape smallest solid superspace (or solid hull)}, of
$A$. We have
\begin{equation}\label{solidcore}
s(A)=\{(b_j)_{j\in \NN}: (b_j \lambda_j)_{j\in \NN}\in A,
\forall\,(\lambda_j)_{j\in \NN}\in l^{\infty}\}
\end{equation} and
\begin{equation}\label{solidalternative}
S(A)=\{(b_j)_{j\in \NN}:\,\exists\,(a_j)_{j\in \NN}\in A:\,|b_j|\le|a_j|,\,\forall\,j\in\NN\},
\end{equation}
e.g. see \cite[p. 594]{BonetTaskinen18}. In our context, the
two following results will show that this notion of solidness is not helping answering the question which $\mathbf{F}$ does belong to the image of $j^{\infty}$ or not (again by identifying a formal weighted power series $\mathbf{F}=\sum_{j=0}^{+\infty}F_jx^j$ by its sequence of coefficients $(F_j)_j$). In particular, we see that the image $j^{\infty}(\mathcal{E}^{0}_{[M]})$ of the Borel map is solid if and only if the Borel map is surjective.


\begin{proposition}\label{solidhull}
Let $M$ be a weight sequence. Then one has $S(j^{\infty}(\mathcal{E}^{0}_{[M]}))=\mathcal{F}_{[M]}$. 
\end{proposition}

\demo{Proof}
Since $j^{\infty}(\mathcal{E}^{0}_{[M]})\subseteq\mathcal{F}_{[M]}$
and $\mathcal{F}_{[M]}$ is solid, we have
$S(j^{\infty}(\mathcal{E}^{0}_{[M]}))\subseteq\mathcal{F}_{[M]}$. For
the proof of the converse inclusion we distinguish between the Roumieu
and the Beurling type.
\medskip

{\itshape Roumieu case.}  Let
$\mathbf{F}=\sum_{j=0}^{+\infty}F_jx^j\in\mathcal{F}_{\{M\}}$ be
given. Then there exist $C,h>0$ such that $|F_j|\le Ch^jm_j$ for all
$j\in\NN$. Let us consider the function
$\theta_{M,C,h}:=C\theta_{(h^jM_j)_j}$ given in Proposition
\ref{prop_thetafunction} using the sequence $(h^jM_j)_j$. By
construction, one has $\mathbf{G}:= \sum_{j=0}^{+
  \infty}\frac{\theta_{M,C,h}^{(j)}(0)}{j!} x^{j}\in
j^{\infty}(\mathcal{E}^{0}_{\{M\}}) $ and
$|\theta_{M,C,h}^{(j)}(0)|/j!\ge C h^{j}m_{j}\geq |F_{j}|$ for every $j \in \NN$.
By \eqref{solidalternative} we have $\mathbf{F}\in
S(j^{\infty}(\mathcal{E}^{0}_{\{M\}}))$ and are done.

\medskip

{\itshape Beurling case.} 
We will apply and recall \cite[Prop. 2.12 (3)]{compositionpaper} (see
also \cite[Proposition 1]{borelmappingquasianalytic} with $\Lambda^1$
instead of $\mathcal{F}$): Since $M$ is a weight sequence, one has
\begin{equation*}\label{roumbeurreduction}
\mathcal{F}_{(M)}=\bigcup_{L\hyperlink{mtriangle}{\vartriangleleft}M, L\in\hyperlink{LCset}{\mathcal{LC}}}\mathcal{F}_{\{L\}},\hspace{30pt}\mathcal{E}^{0}_{(M)}=\bigcup_{L\hyperlink{mtriangle}{\vartriangleleft}M, L\in\hyperlink{LCset}{\mathcal{LC}}}\mathcal{E}^{0}_{\{L\}}.
\end{equation*}
Hence if $\mathbf{F}\in\mathcal{F}_{(M)}$, then
$\mathbf{F}\in\mathcal{F}_{\{L\}}$ for some
$L\in\hyperlink{LCset}{\mathcal{LC}}$ with
$L\hyperlink{mtriangle}{\vartriangleleft}M$. The Roumieu part shows
$\mathbf{F}\in S(j^{\infty}(\mathcal{E}^{0}_{\{L\}}))$ and so, by
$L\hyperlink{mtriangle}{\vartriangleleft}M$, also $\mathbf{F}\in
S(j^{\infty}(\mathcal{E}^{0}_{(M)}))$ follows because $A\subseteq B$
implies $S(A)\subseteq S(B)$. The conclusion follows.  

\qed\enddemo

Concerning the solid core, we have the following result.

\begin{proposition}\label{solidcorelemma}
Let $M$ be a quasianalytic weight sequence such that
$\mathcal{O}^{0}\subsetneq\mathcal{E}^{0}_{[M]}$. Then one has $s(j^{\infty}(\mathcal{E}^{0}_{[M]}))=\mathcal{F}_{\{(j!)_j\}}$.
\end{proposition}

\demo{Proof}
First note that
$\mathcal{F}_{\{(j!)_j\}}=j^\infty(\mathcal{O}^{0})\subseteq
j^\infty(\mathcal{E}^{0}_{[M]})$ and since
$\mathcal{F}_{\{(j!)_j\}}$ is solid by definition, we obtain directly $\mathcal{F}_{\{(j!)_j\}}\subseteq s(j^{\infty}(\mathcal{E}^{0}_{[M]}))$.

\medskip

{\itshape Roumieu case.}  Conversely, let us consider $\mathbf{F}\in
s(j^{\infty}(\mathcal{E}^{0}_{\{M\}}))\subseteq\mathcal{F}_{\{M\}}$
with $\mathbf{F}=\sum_{j=0}^{+\infty}F_jx^j$. Then by
\eqref{solidcore}, one has $|\mathbf{F}|\in
j^{\infty}(\mathcal{E}^{0}_{\{M\}})$  with
$|\mathbf{F}|:=\sum_{j=0}^{+\infty}|F_j|x^j$. Let us assume by
contradiction that $\mathbf{F}\notin\mathcal{F}_{\{(j!)_j\}}$.
We consider $\mathbf{G}\in j^{\infty}(\mathcal{O}^{0})$ with $G_{j}>0$
for every $j \in \NN$, and we set
$$
\mathbf{H} = \sum_{j=0}^{+ \infty} (|F_{j}|+ G_{j} ) x^{j}.
$$
Since $\mathcal{F}_{\{(j!)_j\}}$ is solid, one has
$\mathbf{H}\notin\mathcal{F}_{\{(j!)_j\}}$. Moreover, $\mathbf{H}  \in
\mathcal{F}_{\{M\}}$ and $H_{j}>0$ for every $j \in \NN$. Proposition
\ref{prop_positive} implies that $\mathbf{H}\notin
j^{\infty}(\mathcal{E}^{0}_{\{M\}})$. Using the fact that
$\mathbf{G}\in j^{\infty}(\mathcal{O}^{0}) \subseteq
j^{\infty}(\mathcal{E}^{0}_{\{M\}})$, we obtain
$|\mathbf{F}| = \mathbf{H}- \mathbf{G} \notin
j^{\infty}(\mathcal{E}^{0}_{\{M\}})$, which gives a contradiction.

\medskip

{\itshape Beurling case.} We know that $\mathcal{F}_{\{(j!)_j\}}\subseteq
j^\infty(\mathcal{E}^{0}_{(M)}) \subseteq
j^\infty(\mathcal{E}^{0}_{\{M\}})$, hence $\mathcal{F}_{\{(j!)_j\}}\subseteq
s(j^\infty(\mathcal{E}^{0}_{(M)})) \subseteq
s(j^\infty(\mathcal{E}^{0}_{\{M\}}))$ and the Roumieu case allows to
conclude.

\qed\enddemo

Let us mention that using unions and intersections, the two previous results easily generalize to the case of weight matrices (and so to weight functions by using the associated weight matrix).

\subsection{Characterization of the closedness under the pointwise product}
The aim is now to characterize, as a first step, the closedness of $\mathcal{F}_{[M]}$ and $\mathcal{F}_{[\mathcal{M}]}$ under $\odot$ defined in \eqref{pointwisemulti}. For the weight function case $\mathcal{F}_{[\omega]}$ we need some more preparation and we will study this situation in Section \ref{stabilitypointwiseweightfunction} below in detail.\vspace{6pt}

First we observe that, if $M\in\RR_{>0}^{\NN}$, then one clearly has
that $\mathcal{F}_{[M]}$ is a ring under $\odot$ provided that $M$ has
\begin{equation}\label{pointwise}
\exists\,C,\,h>0\text{ such that }\forall\,j\in\NN\,, \,\,(m_j)^{2}\le Ch^jm_j,
\end{equation}
which is also equivalent to $\sup_{j\in\NN_{>0}}(m_j)^{1/j}<+\infty$
(i.e. $M\hyperlink{mpreceq}{\preceq}(j!)_{j\in\NN}$). 

\medskip

In the general weight matrix setting we consider the following generalizations of \eqref{pointwise}: In the Roumieu case we require
\begin{equation}\label{Roumieusquare}
\forall\,\lambda >0\,\, \exists\,\kappa>0\, , \,\exists\,C,h>0\text{ such
that }\forall\,j\in\NN\, , \,\,(m^{(\lambda)}_j)^2\le Ch^jm^{(\kappa)}_j,
\end{equation}
and in the Beurling case
\begin{equation}\label{Beurlingsquare}
\forall\,\lambda >0\,\, \exists\,\kappa>0\, , \,\exists\,C,h>0\text{ such
that }\forall\,j\in\NN\, ,\,\,(m^{(\kappa)}_j)^2\le Ch^jm^{(\lambda)}_j.
\end{equation}

It is immediate to see that \eqref{Roumieusquare} is preserved under \hyperlink{Mroumapprox}{$\{\approx\}$} and \eqref{Beurlingsquare} under \hyperlink{Mbeurapprox}{$(\approx)$}.

In this situation we can estimate as follows for all $j\in\NN$:
$$|F_jG_j|\le C_1h_1^jm^{(\lambda_1)}_jC_2h_2^jm^{(\lambda_2)}_j\le C_1C_2(h_{1}h_{2})^{j}(m^{(\lambda_3)}_j)^2\le C_1C_2C_3(h_{1}h_2h_{3})^{j}m^{(\kappa)}_j,$$
by taking $\lambda_3:=\max\{\lambda_1,\lambda_2\}$. This shows the
Roumieu case, the Beurling case holds true analogously. So these
conditions are sufficient to have closedness under the pointwise
product. We will show now that under mild additional assumptions on
$\mathcal{M}$, \eqref{Roumieusquare} and \eqref{Beurlingsquare} are
also necessary for the particular case (and thus in the single weight
sequence case \eqref{pointwise}).

\medskip

The proof of the stability of $\mathcal{F}_{\{\mathcal{M}\}}$  under
the pointwise product will use the following classical result, see \cite[Chapitre
I]{mandelbrojtbook} and \cite[Proposition 3.2]{Komatsu73}. Note that
it allows also
to construct the log-convex minorant of a sequence.

\begin{proposition}\label{prop:lcminorant}
Let $M\in\RR_{>0}^{\NN}$ (with $M_0=1$) be a log-convex sequence. Then its
 associated function $\omega_M: \RR_{\ge 0}\rightarrow\RR\cup\{+\infty\}$  defined by
\begin{equation*}\label{assofunc}
\omega_M(t):=\sup_{j\in\NN}\log\left(\frac{t^j}{M_j}\right)\,\,\,\text{
  for }\,t>0,\quad \omega_M(0):=0,
\end{equation*}
satisfies
\begin{equation}\label{lcminorant}
M_j=\sup_{t>0}\frac{t^j}{\exp(\omega_M(t))}\, ,\quad\forall\,j\in\NN.
\end{equation}
\end{proposition}

We say that a family of sequences
$\mathcal{M}=\{M^{(\lambda)}\in\RR_{>0}^{\NN}: \lambda>0\}$ is
\emph{standard log-convex} if $M^{(\lambda)}\le M^{(\kappa)}$ for all
$0<\lambda\le\kappa$ and if
$M^{(\lambda)}\in\hyperlink{LCset}{\mathcal{LC}}$ for all $\lambda>0$
(which is slightly weaker than Definition \ref{defsect41}).

\medskip
We can now state and prove the result of stability under the pointwise
product.

\begin{proposition}\label{pointwisecharacterizationmatrix}
Let $\mathcal{M}=\{M^{(\lambda)}: \lambda>0\}$ be standard log-convex. Then $\mathcal{F}_{\{\mathcal{M}\}}$ is closed under the pointwise product $\odot$ if and only if \eqref{Roumieusquare} holds true and $\mathcal{F}_{(\mathcal{M})}$ is a ring under the product $\odot$ if and only if \eqref{Beurlingsquare} holds true.
\end{proposition}

\demo{Proof}
{\itshape Roumieu case.}  
Assume that $\mathcal{F}_{\{\mathcal{M}\}}$ is  a ring under the
pointwise product and fix an index $\lambda>0$. Since the formal
power series $\mathbf{F}^{\lambda}:=\sum_{j=0}^{+ \infty} m^{(\lambda)}_{j}
  x^j$ belongs to  $\mathcal{F}_{\{\mathcal{M^{(\lambda)}}\}}\subseteq
  \mathcal{F}_{\{\mathcal{M}\}}$, one also has $\mathbf{F}^{\lambda}
  \odot \mathbf{F}^{\lambda} \in \mathcal{F}_{\{\mathcal{M}\}}$. Hence
  there exist an index $\kappa$ and numbers $C,h>0$ such that
$$(m^{(\lambda)}_j)^{2} \le Ch^j m^{(\kappa)}_j$$
for all $j \in \NN$,
and \eqref{Roumieusquare} follows.

\medskip

{\itshape Beurling case.} We follow the ideas from \cite[Section
2]{Bruna} and \cite[Proposition 4.6 $(1)$]{compositionpaper}. We set
$$\mathcal{F}^2_{(\mathcal{M})}:=\left\{\mathbf{F}=\sum_{j=0}^{+ \infty} F_j
x^j\, : \,\,\,\forall\,\lambda>0\,\, \forall\,h>0\,, \,\,(|\mathbf{F}|^2)^{M^{(\lambda)}}_h:=\sup_{j
  \in \NN} \frac{|F_j|^2}{h^{j} m^{(\lambda)}_j} < + \infty\right\}.$$ Note
that both $\mathcal{F}^2_{(\mathcal{M})}$ and
$\mathcal{F}_{(\mathcal{M})}$ are Fr\'echet space spaces under the
canonical projective topology over all $h=h_1^{-1}$ and
$\lambda=\lambda_1^{-1}$, $h_1,\lambda_1\in\NN_{>0}$. By assumption
$\mathcal{F}_{(\mathcal{M})}$ is closed under the pointwise product
which amounts to
$\mathcal{F}_{(\mathcal{M})}\subseteq\mathcal{F}^2_{(\mathcal{M})}$. The
closed graph theorem implies that this last inclusion is
continuous. Consequently, for each $\lambda>0$ and $h>0$, there exist
$\kappa>0$ and $C,h_1>0$ such that for each
$\mathbf{F}=\sum_{j=0}^{+\infty}F_jx^j\in\mathcal{F}_{(\mathcal{M})}$,
\begin{equation}\label{eq:continuity}
(|\mathbf{F}|^2)^{M^{(\lambda)}}_{h}=\sup_{j\in\NN}\frac{|F_j|^2}{h^j
    m^{(\lambda)}_j}\le
  C\sup_{j\in\NN}\frac{|F_j|}{h_1^jm^{(\kappa)}_j}=C|\mathbf{F}|^{M^{(\kappa)}}_{h_1}.
\end{equation}

For every $s\ge 0$, let us consider the function
$f_s(t):=\sin(st)+\cos(st)$, $t\in\RR$, and let us show  that
$\mathbf{F}^s:=\sum_{j=0}^{+\infty}\frac{f_s^{(j)}(0)}{j!}x^j\in\mathcal{F}_{(\mathcal{M})}$. Indeed,
if $s>0$ (the case $s=0$ is obvious),
note that $|f_s^{(j)}(0)|=s^j$ for all $j\in\NN$ and since for all
$\lambda>0$, $(M^{(\lambda)}_j)^{1/j}\rightarrow+\infty$ as $j\rightarrow+\infty$, it is direct to check that for all  $\lambda>0$ and
all $h>0$ there exists some $C>0$ such that $ s^j\le
C h^jM^{(\lambda)}_j$ for all $j \in \NN$. 
Now, inequality \eqref{eq:continuity} applied to the family $\mathbf{F}^s$, $s\ge 0$, and with the choice $h=1$ yields
\begin{align*}
\sup_{j\in\NN}\frac{s^{2j}}{\widehat{M}^{(\lambda)}_j}&=\sup_{j\in\NN}\frac{|f_s^{(j)}(0)|^2}{j!^2m^{(\lambda)}_j}=\sup_{j\in\NN}\frac{|F^s_j|^2}{m^{(\lambda)}_j}\le C\sup_{j\in\NN}\frac{|F^s_j|}{h_1^jm^{(\kappa)}_j}=C\sup_{j\in\NN}\frac{|f_s^{(j)}(0)|}{j!h_1^j m^{(\kappa)}_j}=C\sup_{j\in\NN}\frac{s^{j}}{h_1^jM^{(\kappa)}_j},
\end{align*}
where we have put
$\widehat{M}^{(\lambda)}:=(j!M^{(\lambda)}_j)_{j\in\NN}$. This
implies in turn $\exp(\omega_{\widehat{M}^{(\lambda)}}(s^2))\le
C\exp(\omega_{M^{(\kappa)}}(s/h_1))$ for all $s\ge 0$,  where the
associated function is defined in Proposition \ref{prop:lcminorant}. Using
\eqref{lcminorant}
we get for all $j\in\NN$:
\begin{align*}
\widehat{M}^{(\lambda)}_{j}=\sup_{t\ge
  0}\frac{t^{j}}{\exp(\omega_{\widehat{M}^{(\lambda)}}(t))} &
 =\sup_{t\ge  0}\frac{t^{2j}}{\exp(\omega_{\widehat{M}^{(\lambda)}}(t^2))}\\
& \ge\frac{1}{C}\sup_{t\ge
  0}\frac{t^{2j}}{\exp(\omega_{M^{(\kappa)}}(t/h_1))}\\
& =\frac{h_1^{2j}}{C}\sup_{t\ge 0}\frac{t^{2j}}{\exp(\omega_{M^{(\kappa)}}(t))}
 =\frac{h_1^{2j}}{C}M^{(\kappa)}_{2j}.
\end{align*}
Consequently $M^{(\kappa)}_{2j}\le C
h_1^{-2j}\widehat{M}^{(\lambda)}_j$ for all $j\in\NN$ follows. Using
the log-convexity of $M^{(\kappa)}$, one knows that the sequence
$(M^{(\kappa)}_{j})^{1/j})_{j}$ is increasing, hence
$(M^{(\kappa)}_j)^2\le M^{(\kappa)}_{2j}$ for all $j \in \NN$. This finally yields $j!^2(m^{(\kappa)}_j)^2=(M^{(\kappa)}_j)^2\le M^{(\kappa)}_{2j}\le C h_1^{-2j}\widehat{M}^{(\lambda)}_j=Ch_1^{-2j}j!^2m^{(\lambda)}_j$ and so \eqref{Beurlingsquare} follows.
\qed\enddemo

\begin{remark}\label{nonsenseweightsequence}
Consequently, if $\mathcal{M}$ is standard log-convex and constant and so we deal with $M\in\hyperlink{LCset}{\mathcal{LC}}$, then
$\mathcal{F}_{\{M\}}$ and/or $\mathcal{F}_{(M)}$ is a ring under the
pointwise product if and only if
$\sup_{j\in\NN_{>0}}(m_j)^{1/j}<+\infty$ which precisely means
$\mathcal{E}_{\{M\}}\subseteq\mathcal{C}^{\omega}$
resp. $\mathcal{E}_{(M)}\subseteq\mathcal{C}^{\omega}$ (e.g. see
\cite[Proposition 4.6]{compositionpaper}). But this is a situation
which cannot be considered under the assumptions of the main result Theorem \ref{thm:algebraabinfinity} of Section \ref{infinityalgebraweightsequence} above (\eqref{strictincl} is violated). Note that $\sup_{j\in\NN_{>0}}(m_j)^{1/j}<+\infty$ is clearly stable under $\hyperlink{approx}{\approx}$ and if $M$ is a weight sequence in the sense of Definition \ref{def_weightseq}, then $\mathcal{F}_{\{M\}}$ and/or $\mathcal{F}_{(M)}$ is a ring under the
pointwise product if and only if $M\hyperlink{approx}{\approx}(j!)_j$ (by combining $(III)$ and $\sup_{j\in\NN_{>0}}(m_j)^{1/j}<+\infty$) and so if and only if $\mathcal{E}_{\{M\}}=\mathcal{C}^{\omega}$.
\end{remark}

Instead of \eqref{Roumieusquare}
resp. \eqref{Beurlingsquare}, it would have been natural to assume on  $\mathcal{M}$ also the following assumptions:
\begin{equation}\label{Roumieubigsquare}
\forall\,\lambda >0\,\, \exists\,\kappa>0\, , \,\exists\,C,h>0\text{ such
that }\forall\,j\in\NN\, , \,\,(M^{(\lambda)}_j)^2\le Ch^jM^{(\kappa)}_j,
\end{equation}
resp.
\begin{equation}\label{Beurlingbigsquare}
\forall\,\lambda >0\,\, \exists\,\kappa>0\, , \,\exists\,C,h>0\text{ such
that }\forall\,j\in\NN\, , \,\,(M^{(\kappa)}_j)^2\le Ch^jM^{(\lambda)}_j.
\end{equation}
\eqref{Roumieubigsquare} is preserved under \hyperlink{Mroumapprox}{$\{\approx\}$} and \eqref{Beurlingbigsquare} under \hyperlink{Mbeurapprox}{$(\approx)$}.

Note that $\eqref{Roumieubigsquare}\Rightarrow\eqref{Roumieusquare}$
resp. $\eqref{Beurlingbigsquare}\Rightarrow\eqref{Beurlingsquare}$
whereas the equivalences will fail in general, see also the example in
Section \ref{counterexample} below.

\subsection{Example of a quasianalytic weight matrix}\label{counterexample}
In contrast to the single weight sequence case we will construct now
an example which shows that \eqref{Roumieusquare} and/or
\eqref{Beurlingsquare} can even hold true for quasianalytic weight
matrices $\mathcal{M}$ satisfying
$\mathcal{C}^{\omega}\subsetneq\mathcal{E}_{[\mathcal{M}]}$, i.e. for
$\mathcal{M}$ having \eqref{strictincl2}. So this weight matrix satisfies the requirements of Theorem
\ref{thm:algebraabinfinityfullmatrix} and hence it illustrates that in
the general matrix setting an equivalent of Theorem
\ref{thm:algebraabinfinityfullmatrix} using the pointwise product
makes sense, see Theorem
\ref{thm:algebraabinfinitymatrixpointwise} below.

\medskip

For this we consider the matrix $\mathcal{M}:=\{M^{(\lambda)}\in
\RR^{{\NN}}_{>0}: \lambda>0\}$ with each $M^{(\lambda)}$ defined by
its quotients
$\mu^{(\lambda)}_j:=\frac{M^{(\lambda)}_j}{M^{(\lambda)}_{j-1}}$ as
follows: Let $j_0\in\NN$ be the smallest integer satisfying
$\log(\log(j))\ge 1$ for all $j>j_0$ (and so not depending on
$\lambda$) and
put $$1=\mu^{(\lambda)}_0=\dots=\mu^{(\lambda)}_{j_0},\hspace{20pt}\mu^{(\lambda)}_j=j\big(\log(\log(j))\big)^{\lambda},\,\forall
 j>j_0.$$
So $j\mapsto\mu^{(\lambda)}_j$ is increasing for each $\lambda>0$,
i.e. each $M^{(\lambda)}$ is log-convex, and even
$\lim_{j\rightarrow+\infty}\mu^{(\lambda)}_j/j=+\infty$ for each
$\lambda>0$ is valid. It is known that this also implies
$\lim_{j\rightarrow+\infty}(m^{(\lambda)}_j)^{1/j}=+\infty$ for each
$\lambda>0$ (e.g. see the argument given on
\cite[p. 104]{compositionpaper}), hence $\mathcal{M}$ is a weight
matrix and satisfies both requirements in
\eqref{strictincl2} (and consequently \eqref{pointwise} does not hold true for any $M^{(\lambda)}$). Moreover, $\mathcal{M}$ is quasianalytic because
each $M^{(\lambda)}$ is clearly quasianalytic.

\medskip

Let us now show that both \eqref{Roumieusquare} and \eqref{Beurlingsquare}
hold true. For all $j>j_0$, one has
$$(m^{(\lambda)}_j)^2=\prod_{i=1}^j\left(\frac{\mu^{(\lambda)}_i}{i}\right)^2=\prod_{i=1}^{j_0}\frac{1}{i^2}\prod_{i=j_0+1}^j\big(\log(\log(i))\big)^{2\lambda}\le \prod_{i=1}^{j_0}\frac{1}{i}\prod_{i=j_0+1}^j\big(\log(\log(i))\big)^{\kappa}=m^{(\kappa)}_j,$$
by taking $\kappa:=2\lambda$ resp. $\lambda:=\kappa/2$.

It is also immediate to see $\lim_{j\rightarrow+\infty}\frac{\mu_j^{(\kappa)}}{\mu_j^{(\lambda)}}=+\infty$ for all $0<\lambda<\kappa$ which implies that all sequences are pairwise not equivalent because $M^{(\lambda)}\hyperlink{mtriangle}{\vartriangleleft}M^{(\kappa)}$ for all $0<\lambda<\kappa$.
\medskip


\vspace{6pt}

Note that $\mathcal{M}$ violates both \eqref{Roumieubigsquare} and
\eqref{Beurlingbigsquare}. Indeed, for all $j>j_0$ we have
\begin{eqnarray*}
&&(M^{(\lambda)}_j)^2\le Ch^jM^{(\kappa)}_j\\
& \Leftrightarrow & \prod_{i=j_0+1}^ji^2\big(\log(\log(i))\big)^{2\lambda}\le C h^j\prod_{i=j_0+1}^ji\big(\log(\log(i))\big)^{\kappa}
\\
& \Leftrightarrow & \prod_{i=j_0+1}^ji\le
                  Ch^j\prod_{i=j_0+1}^j\big(\log(\log(i))\big)^{\kappa-2\lambda}\\
& \Leftrightarrow &  j!\le j_0!Ch^j\prod_{i=j_0+1}^j\big(\log(\log(i))\big)^{\kappa-2\lambda}.
\end{eqnarray*}
But this cannot hold true for all $j\in\NN$ for any given numbers $C$
and $h$ large, since, by Stirling's formula, the left-hand side is
increasing like $j\mapsto \left(\frac{j}{e}\right)^{j}\sqrt{2\pi j}$, whereas the right-hand side is bounded by above by $j_0!Ch^j\log(\log(j))^{j(\kappa-2\lambda)}$.


\vspace{6pt}
It shall be noted that, by the characterization shown in Proposition \ref{pointwisecharacterizationmatrix}, we have stability under $\odot$ for both $\mathcal{F}_{\{\mathcal{M}\}}$ and $\mathcal{F}_{(\mathcal{M})}$. However, even in this situation it is still impossible to obtain closedness under $\odot$ for $j^{\infty}(\mathcal{E}^{0}_{\{\mathcal{M}\}})$: Take $\theta_{M^{(\lambda_0)}}$ for some $\lambda_0>0$ and put $\mathbf{F}:=j^{\infty}(\theta_{M^{(\lambda_0)}})$. Then clearly $\mathbf{F}\in\mathcal{F}_{\{M^{(\lambda_0)}\}}\subseteq\mathcal{F}_{\{\mathcal{M}\}}$ but $|\mathbf{F}|\notin j^{\infty}(\mathcal{E}^{0}_{\{\mathcal{M}\}})$ (with $|\mathbf{F}|:=\sum_{j=0}^{+\infty}|F_j|x^j$) since $|\mathbf{F}|\notin j^{\infty}(\mathcal{E}^{0}_{\{L\}})$ for any quasianalytic weight sequence $L$ (see Proposition \ref{prop_positive}) and so in particular this holds true for the sequence $L$ coming from Lemma \ref{roumieuauxiliary}.\vspace{6pt}

We close this section with the following observation: Not for all (quasianalytic) weight matrices the characterizing conditions \eqref{Roumieusquare} and \eqref{Beurlingsquare} are satisfied simultaneously.

For this we consider $\mathcal{N}:=\{(j!)_j, M^{(\lambda_0)}\}$ with $M^{(\lambda_0)}$ denoting one of the sequences belonging to the matrix $\mathcal{M}$ constructed above. So $\mathcal{N}$ is a weight matrix consisting only of two non-equivalent (quasianalytic) weight sequences and so $\mathcal{F}_{(\mathcal{N})}=\mathcal{F}_{((j!)_j)}$, $\mathcal{F}_{\{\mathcal{N}\}}=\mathcal{F}_{\{M^{(\lambda_0)}\}}$. Then \eqref{Beurlingsquare}, which amounts to \eqref{pointwise} for $(j!)_j$ holds true, whereas \eqref{Roumieusquare} for $\mathcal{N}$, i.e. \eqref{pointwise} for $M^{(\lambda_0)}$, fails. Note that $j!\le M^{(\lambda_0)}_j$ only holds true for all $j\in\NN$ large, but $M^{(\lambda_0)}$ can be replaced by an equivalent sequence satisfying this pointwise estimate for all $j\in\NN$ (as required in Definition \ref{defsect41}) and defining the same matrix.

\subsection{Algebrability for the general matrix setting}\label{weightmatrixpointwise}

As seen by the example constructed in Section \ref{counterexample}, in
the general weight matrix setting it makes also sense to consider on
$\mathcal{F}_{[\mathcal{M}]}$ the pointwise product. We show the
following result analogous to Theorem \ref{thm:algebraabinfinityfullmatrix} for the convolution product but the proof will simplify at several steps due to the fact that multiplying two lacunary series w.r.t. $\odot$ does not change and mix the indices $j\in\NN$ with $F_j\neq 0$.

\begin{theorem}\label{thm:algebraabinfinitymatrixpointwise}
Let $\mathcal{M}, \mathcal{N}$ be two quasianalytic weight
matrices.
We assume
\begin{itemize}
\item[$(i)$] in the Roumieu case that $\mathcal{N}$ satisfies
  \eqref{Roumieusquare} and
  $\mathcal{O}^{0}\subsetneq\mathcal{E}^{0}_{\{\mathcal{N}\}}$,
\item[$(ii)$] in the Beurling case that $\mathcal{N}$ satisfies \eqref{Beurlingsquare} and $\mathcal{O}^{0}\subsetneq\mathcal{E}^{0}_{(\mathcal{N})}$.
\end{itemize}
Then $\mathcal{F}_{[\mathcal{N}]} \setminus
j^{\infty}(\mathcal{E}^{0}_{\{\mathcal{M}\}})$ is $\mathfrak{c}$-algebrable in
$\mathcal{F}_{[\mathcal{N}]}$ endowed with the pointwise product (hence $\mathcal{F}_{[\mathcal{N}]} \setminus j^{\infty}(\mathcal{E}^{0}_{(\mathcal{M})})$ too).
\end{theorem}

\demo{Proof}
As in the proof of  Theorem \ref{thm:algebraabinfinityfullmatrix}, one
can use Lemma \ref{roumieuauxiliary} to reduce the proof to the case
of a  quasianalytic weight sequence $L$ instead of $\mathcal{M}$.
By assumption, one can construct an increasing sequence $(k_p)_{p\in \NN} $ of natural numbers satisfying
\begin{enumerate}[(i)]
\item $k_0=1$ and $k_p > k_{p-1}$ for every $p \in \NN_{>0}$,
\item $\lim_{p\rightarrow + \infty}\left(n^{(1/(p+1))}_{k_p}\right)^\frac{1}{k_p}=+\infty$,
\item $\sum_{j=0}^{k_{p-1}}\left|\omega^L_{j,k_p}-1\right|n_j^{(p)}\le 1$ for every $p \in \NN_{>0}$.
\end{enumerate}
We proceed then exactly as in the proof of Theorem
\ref{thm:algebraabinfinityfullmatrix} to construct formal power series
$\mathbf{F}^{b}$, $b \in \mathcal{H}$, and we remark that if
$$\mathbf{G} = \sum_{l=1}^{L'}\alpha_l \underbrace{(\mathbf{F}^{b_1}\odot
  \dots \odot\mathbf{F}^{b_1})}_{i_{l,1} \text{ times }} \odot \dots
\odot \underbrace{(\mathbf{F}^{b_J}\odot
  \dots \odot\mathbf{F}^{b_J})}_{i_{l,J} \text{ times }},$$
then
$$
G_{j } =
\begin{cases}
\displaystyle\sum_{l=1}^{L'} \alpha_l
\big(n_{k_p}^{(1/(p+1))}\big)^{{i_{l,1}}{b_1} + \cdots  +
  {i_{l,J}}{b_J}} & \text{ if } j= k_{p}\\[2ex]
0  & \text{ if } j\notin\{k_{p}:p \in \NN\}.
\end{cases}
$$
To conclude, one follows the same ideas as in the proofs of Theorem
\ref{thm:algebraabinfinity} and Theorem
\ref{thm:algebraabinfinityfullmatrix}.

\qed\enddemo

The identity for $\odot$ is given by $\mathbf{E}_{\odot}=\sum_{j=0}^{+\infty}1x^j$ and so $\mathbf{E}=j^{\infty}(f)$ with $f(x):=\sum_{j=0}^{+\infty}x^j$ representing a real analytic germ at $0$. Consequently also in this setting each algebra contained in $\mathcal{F}_{[\mathcal{N}]} \setminus j^{\infty}(\mathcal{E}^{0}_{\{\mathcal{M}\}})$ does not contain the identity $\mathbf{E}_{\odot}$ anymore.

\section{On the stability under the pointwise product of $\mathcal{F}_{[\omega]}$}\label{stabilitypointwiseweightfunction}
The goal of this Section is to show that, similarly as commented in
Remark \ref{nonsenseweightsequence} for the single weight sequence
situation, the problem of algebrability with respect to $\odot$ cannot be considered for $\mathcal{F}_{[\omega]}$ within the quasianalytic setting. More precisely we will show that all required assumptions on $\omega$ can never be satisfied simultaneously. While in the
weight function case we can have the situation that $\mathcal{F}_{[\Omega]}=\mathcal{F}_{[\omega]}$ is closed under the
pointwise product $\odot$ and $\mathcal{E}_{[\Omega]}=\mathcal{E}_{[\omega]}$
is strictly containing the real analytic functions, we will see below that this situation forces already non-quasianalyticity for $\omega$. Consequently the matrix constructed in Section \ref{counterexample} above cannot be associated with a weight function $\omega$.




\vspace{6pt}

In order to do so first recall that, as shown in Lemma \ref{pointwisecharacterizationmatrix} above, \eqref{Roumieusquare} resp. \eqref{Beurlingsquare} are characterizing the closednees under the pointwise product for $\mathcal{F}_{\{\Omega\}}=\mathcal{F}_{\{\omega\}}$ resp. $\mathcal{F}_{(\Omega)}=\mathcal{F}_{(\omega)}$. Hence we have to show which condition on $\omega$ guarantees that $\Omega$ satisfies \eqref{Roumieusquare} resp. \eqref{Beurlingsquare} and for this we have to introduce some notation and recall several results.\vspace{6pt}

Let $\omega$ be given satisfying \hyperlink{om0}{$(\omega_0)$}, \hyperlink{om3}{$(\omega_3)$} and \hyperlink{om4}{$(\omega_4)$}, then as shown in \cite[Section 5]{compositionpaper}, respectively \cite[Theorem 4.0.3, Lemma 5.1.3]{dissertation} and reproved in \cite[Lemma 2.5]{sectorialextensions} in a more precise way, we have
\begin{equation}\label{goodequivalence}
\forall\,\lambda>0\,\,\exists\,C_{\lambda}>0\text{ such that
}\forall\,t\ge 0\, , \,\,\lambda\omega_{W^{(\lambda)}}(t)\le\omega(t)\le 2\lambda\omega_{W^{(\lambda)}}(t)+C_{\lambda}.
\end{equation}
In particular we have $\omega\hyperlink{sim}{\sim}\omega_{W^{(\lambda)}}$ for all $\lambda>0$.\vspace{6pt}

Moreover, for any $h:(0,+\infty)\rightarrow[0,+\infty)$ which is nonincreasing and such that $\lim_{s\rightarrow 0}h(s)=+\infty$, we can define the so-called {\itshape lower Legendre conjugate (or envelope)} $h_{\star}:[0,+\infty)\rightarrow[0,+\infty)$ of $h$ by
\begin{equation*}\label{omegaconjugate0}
h_{\star}(t):=\inf_{s>0}\{h(s)+ts\},\hspace{15pt}t\ge 0.
\end{equation*}
We are summarizing some facts for this conjugate, see also
\cite[Section 3.1]{sectorialextensions}. The function $h_{\star}$ is
clearly nondecreasing, continuous and concave, and
$\lim_{t\rightarrow+\infty}h_{\star}(t)=+\infty$, see \cite[(8),
p. 156]{Beurling72}. Moreover, if $\lim_{s\to +\infty}h(s)=0$ then
$h_{\star}(0)=0$, and so $h_{\star}$ satisfies all properties from
\hyperlink{om0}{$(\omega_0)$} except normalization. In the forthcoming proof we will apply this conjugate to $h(t)=\omega^{\iota}(t):=\omega(1/t)$, where $\omega$ is a weight function, so that $(\omega^{\iota})_{\star}$ is again a weight function (except normalization); in particular, we will frequently find the case $h(t)=\omega^{\iota}_M(t)=\omega_M(1/t)$ for $M\in\RR_{>0}^{\NN}$ with $\lim_{p\rightarrow+\infty}(M_p)^{1/p}=+\infty$.

\medskip

Now we are able to formulate the first main characterizing result.

\begin{theorem}\label{mixedomega7}
Let $\omega$ be given satisfying \hyperlink{om0}{$(\omega_0)$},
\hyperlink{om3}{$(\omega_3)$} and \hyperlink{om4}{$(\omega_4)$}, and let $\Omega=\{W^{(\lambda)}: \lambda>0\}$ be the matrix associated with $\omega$. Then $\Omega$ satisfies \eqref{Roumieusquare} and/or \eqref{Beurlingsquare} if and only if
\begin{equation}\label{mixedomega7equ}
\exists\,H>0\,\exists\,C>0\text{ such that }\forall\,t\ge 0\, , \,\,(\omega^{\iota})_{\star}(t^2)\le C\omega(H t)+C.
\end{equation}
\end{theorem}

Consequently, if $\omega\in\hyperlink{omset1}{\mathcal{W}}$, then \eqref{mixedomega7equ} is equivalent to having that $\mathcal{F}_{\{\omega\}}=\mathcal{F}_{\{\Omega\}}$ and/or $\mathcal{F}_{(\omega)}=\mathcal{F}_{(\Omega)}$ is closed under the pointwise product $\odot$.

\demo{Proof}
First, let us assume that $\Omega$ satisfies \eqref{Roumieusquare}
and/or \eqref{Beurlingsquare} with indices $\lambda$ and $\kappa$. We
will prove here the Roumieu case, the Beurling case can be treated in
a similar way. If we put
$\widehat{W}^{(\lambda)}:=(j!W^{(\lambda)}_j)_{j\in\NN}$, then we have
 $(W^{(\lambda)}_j)^2\le
 Ch^jj!^2w^{(\kappa)}_j=Ch^j\widehat{W}^{(\kappa)}_j$. Hence for all $t\ge 0$ and $j\in\NN$ we get $\frac{t^j}{\widehat{W}^{(\kappa)}_j}\le C\frac{(ht)^j}{(W^{(\lambda)}_j)^2}=C\left(\frac{(\sqrt{ht})^j}{W^{(\lambda)}_j}\right)^2$, and applying logarithm to this inequality yields $\omega_{\widehat{W}^{(\kappa)}}(t)\le 2\omega_{W^{(\lambda)}}(\sqrt{ht})+\log(C)$.\vspace{6pt}

From \cite[Lemma 3.4 $(ii)$, $(3.6)$]{sectorialextensions} applied to $Q=M=\widehat{W}^{(\kappa)}$ (recall that $W^{(\kappa)}\in\Omega$), we know that
\begin{equation}\label{Dynkinequiv3}
\forall\,t\ge\frac{\widehat{W}^{(\kappa)}_1}{\widehat{W}^{(\kappa)}_0}\,
,
\,\,\omega_{\widehat{W}^{(\kappa)}}(t)\le(\omega^{\iota}_{W^{(\kappa)}})_{\star}(t)\le 1+\omega_{\widehat{W}^{(\kappa)}}(et).
\end{equation}
The second inequality of \eqref{Dynkinequiv3} yields
$$(\omega^{\iota}_{W^{(\kappa)}})_{\star}(t)\le 1+\omega_{\widehat{W}^{(\kappa)}}(et)\le 1+2\omega_{W^{(\lambda)}}(\sqrt{het})+\log(C).$$
By using the first inequality of \eqref{goodequivalence} we see for all $t\ge 0$ that $2\omega_{W^{(\lambda)}}(\sqrt{het})\le\frac{2}{\lambda}\omega(\sqrt{het})$ and the second inequality of \eqref{goodequivalence} implies $(\omega^{\iota}_{W^{(\kappa)}})_{\star}(t)=\inf_{s>0}\{\omega^{\iota}_{W^{(\kappa)}}(s)+st\}\ge\inf_{s>0}\{\frac{1}{2\kappa}\omega^{\iota}(s)+st\}-\frac{C_{\kappa}}{2\kappa}=\frac{1}{2\kappa}(\omega^{\iota})_{\star}(2\kappa t)-\frac{C_{\kappa}}{2\kappa}$. Thus, combining everything, we have shown for all $t$ (large enough) that
$$(\omega^{\iota})_{\star}(t^2)\le\frac{4\kappa}{\lambda}\omega(\sqrt{he/(2\kappa)}t)+2\kappa(1+\log(C))+C_{\kappa},$$
hence \eqref{mixedomega7equ} is satisfied.

\medskip

Conversely, assume now that \eqref{mixedomega7equ} holds true with constants $C>0$ and $H>0$. First, let in the following computations $\lambda, \kappa>0$ be arbitrary but fixed. The second inequality of \eqref{goodequivalence} yields $C\omega(Ht)+C\le 2C\lambda\omega_{W^{(\lambda)}}(Ht)+C(C_{\lambda}+1)$, whereas the first one implies $(\omega^{\iota})_{\star}(t^2)=\inf_{s>0}\{\omega^{\iota}(s)+st^2\}\ge\inf_{s>0}\{\kappa\omega^{\iota}_{W^{(\kappa)}}(s)+st^2\}=\kappa(\omega^{\iota}_{W^{(\kappa)}})_{\star}(t^2/\kappa)$. Moreover, the first estimate in \eqref{Dynkinequiv3} implies $\kappa(\omega^{\iota}_{W^{(\kappa)}})_{\star}(t^2/\kappa)\ge\kappa\omega_{\widehat{W}^{(\kappa)}}(t^2/\kappa)$ for all $t\ge\frac{\widehat{W}^{(\kappa)}_1}{\widehat{W}^{(\kappa)}_0}$, and altogether
$$\exists\,D=D_{C,\lambda,\kappa}\text{ such that }\forall\,t\ge 0\, ,
\,\,\omega_{\widehat{W}^{(\kappa)}}(t^2)\le\frac{2C\lambda}{\kappa}\omega_{W^{(\lambda)}}(H\sqrt{\kappa}t)+D.$$
Now take $\kappa=C\lambda$ and with this choice, by using Proposition
\ref{prop:lcminorant}, we can estimate as follows for all $j\in\NN$
\begin{align*}
\widehat{W}^{(\kappa)}_j=\sup_{t\ge
  0}\frac{t^j}{\exp(\omega_{\widehat{W}^{(\kappa)}}(t))}
& =\sup_{t\ge
  0}\frac{t^{2j}}{\exp(\omega_{\widehat{W}^{(\kappa)}}(t^2))}\\
& \ge\frac{1}{\exp(D)}\sup_{t\ge 0}\frac{t^{2j}}{\exp(2\omega_{W^{(\lambda)}}(H\sqrt{\kappa}t))}
\\
& =\frac{1}{\exp(D)H^{2j}\kappa^j}\left(\sup_{t\ge 0}\frac{t^j}{\exp(\omega_{W^{(\lambda)}}(t))}\right)^2=\frac{1}{\exp(D)H^{2j}\kappa^j}(W^{(\lambda)}_j)^2,
\end{align*}
hence $(W^{(\lambda)}_j)^2\le\exp(D)(H^2\kappa)^jj!W^{(\kappa)}_j$ for
all $j\in\NN$. This proves both \eqref{Roumieusquare} and
\eqref{Beurlingsquare} since $C\lambda=\kappa$ and $C$ is only
depending on given $\omega$.
\qed\enddemo

The characterizing property \eqref{mixedomega7equ} is looking similar to the following growth property on $\omega$, see \cite{dissertation}, \cite{Franken95}, \cite[Theorem 5.14 $(4)$]{compositionpaper}
(called $(\omega_8)$ in there) and \cite[Appendix A]{sectorialextensions1} (denoted by $(\omega_7)$ there):
\begin{equation}\label{omega7}
\exists\,H>0\,\,\exists\,C>0\text{ such that }\forall\,t\ge 0\, , \,\,\omega(t^2)\le C\omega(H t)+C.
\end{equation}
For any $\omega$ satisfying \hyperlink{om0}{$(\omega_0)$}, \hyperlink{om3}{$(\omega_3)$} and \hyperlink{om4}{$(\omega_4)$} condition \eqref{omega7} does always imply \hyperlink{om1}{$(\omega_1)$}, see \cite[Appendix A]{sectorialextensions1}.

In \cite[Lemma A.1]{sectorialextensions1} it has been shown that for any
$\omega\in\hyperlink{omset1}{\mathcal{W}}$ with \eqref{omega7} the associated matrix $\Omega$ does have both \eqref{Roumieubigsquare} and
\eqref{Beurlingbigsquare} (by having a precise relation between the indices $\lambda$ and $\kappa$). Following the proof of \cite[Lemma A.1
$(ii)\Rightarrow(i)$]{sectorialextensions1} and replacing $Al$ by $l_1$ there it is straightforward to see that \eqref{Roumieubigsquare} and/or
\eqref{Beurlingbigsquare} are implying \eqref{omega7}, see \cite[Lemma 5.4.1]{dissertation} and also the first half of the proof of Theorem \ref{mixedomega7} (in fact for this implication one only needs that the inequalities in \eqref{Roumieubigsquare} or
\eqref{Beurlingbigsquare} are valid for some pair of indices $\lambda$ and $\kappa$).

Thus any matrix $\Omega$ associated with a function $\omega$ satisfying \eqref{omega7} will always have both \eqref{Roumieusquare} and
\eqref{Beurlingsquare}, too.\vspace{6pt}

However, \eqref{omega7} implies quite strong, and in our situation undesired, properties for the
associated weight matrix $\Omega$. More precisely, by the results shown in \cite[Appendix A]{sectorialextensions1} we have that for any $\omega$ satisfying \hyperlink{om0}{$(\omega_0)$}, \hyperlink{om3}{$(\omega_3)$} and \hyperlink{om4}{$(\omega_4)$} property \eqref{omega7} does imply the {\itshape strong non-quasianalyticity condition} for weight functions
\begin{equation}\label{omegasnq}
\exists\,C>0\text{ such that }\forall\,y>0\, , \,\, \int_1^{+\infty}\frac{\omega(y t)}{t^2}dt\le C\omega(y)+C,
\end{equation}
and so in particular \eqref{omegaQ} has to fail.

By the results shown in \cite{BonetBraunMeiseTaylorWhitneyextension} (see also \cite{BonetMeiseTaylorSurjectivity}) it follows that for given $\omega\in\hyperlink{omset1}{\mathcal{W}}$ condition \eqref{omegasnq} is characterizing $j^{\infty}(\mathcal{E}^{0}_{[\omega]})=\mathcal{F}_{[\omega]}$, i.e. the surjectivity of the Borel mapping. Note that in \cite{BonetBraunMeiseTaylorWhitneyextension} and \cite{BonetMeiseTaylorSurjectivity} non-quasianalyticity for $\omega$ was a basic assumption but which is superfluous provided that $\mathcal{O}^0\subsetneq\mathcal{E}^{0}_{[\omega]}$ (characterized by \eqref{strictinclusioncomega}): On the one hand it is clear that \eqref{omegasnq} forces non-quasianalyticity for $\omega$. On the other hand, if $\omega\in\hyperlink{omset1}{\mathcal{W}}$ with $\mathcal{O}^0\subsetneq\mathcal{E}^{0}_{[\omega]}$ and $j^{\infty}(\mathcal{E}^{0}_{[\omega]})=\mathcal{F}_{[\omega]}$ then $\omega$ has to be non-quasianalytic: If $\omega$ would be quasianalytic, this would contradict \cite[Cor. 2, Cor. 4]{borelmappingquasianalytic} (similarly see also \cite[Cor. 2]{whitneyextensionmixedweightfunctionII}).\vspace{6pt}

We are gathering now some more observations.\vspace{6pt}

\begin{itemize}
\item[$(i)$] Under the assumptions of Theorem \ref{mixedomega7}, one has that \eqref{mixedomega7equ} and
$\omega\hyperlink{sim}{\sim}(\omega^{\iota})_{\star}$ hold true if and
only if \eqref{omega7} holds true. Indeed, \eqref{mixedomega7equ} together with
$\omega\hyperlink{sim}{\sim}(\omega^{\iota})_{\star}$ immediately
imply \eqref{omega7}. For the converse, first note that \cite[Lemma 5.1, Corollary 5.2]{sectorialextensions} can be applied to each $\omega$ as assumed in the result with \eqref{omega7} since \hyperlink{om1}{$(\omega_1)$} follows as mentioned above. Hence we get $\omega_{\widehat{W}^{(\lambda)}}\hyperlink{sim}{\sim}(\omega^{\iota})_{\star}\hyperlink{sim}{\sim}\omega_{\widehat{W}^{(\kappa)}}$ for all $\lambda,\kappa>0$.
By \cite[Lemma A.1]{sectorialextensions} we know that
$$\forall\,\lambda>0\,\exists\,\kappa>0\,\exists\,C\ge 1\text{ such
  that }\forall\,j\in\NN\, , \,\,\widehat{W}^{(\lambda)}_j=j!W^{(\lambda)}_j\le C^jW^{(\kappa)}_j,$$
hence $\omega_{W^{(\kappa)}}(t)\le\omega_{\widehat{W}^{(\lambda)}}(Ct)\le\omega_{W^{(\lambda)}}(Ct)$. By \eqref{goodequivalence} and the fact that $\omega$ has \hyperlink{om1}{$(\omega_1)$} we have shown $\omega\hyperlink{sim}{\sim}(\omega^{\iota})_{\star}$. Obviously this and \eqref{omega7} together imply \eqref{mixedomega7equ}.
\item[$(ii)$] Let $\omega\in\hyperlink{omset1}{\mathcal{W}}$ be given. Then $\omega\hyperlink{sim}{\sim}(\omega^{\iota})_{\star}$ implies $\gamma(\omega)=+\infty$, with $\gamma$ denoting the growth index studied in detail in \cite{firstindexpaper} and used in the extension results in \cite{sectorialextensions}, \cite{sectorialextensions1} (the fact that $\omega$ has \hyperlink{om1}{$(\omega_1)$} is equivalent to having $\gamma(\omega)>0$, see \cite[Corollary 2.14]{firstindexpaper}). To show this note that by \cite[Proposition 2.22, Corollary 2.26]{firstindexpaper} we have $\gamma(\omega)+1=\gamma((\omega^{\iota})_{\star})=\gamma(\omega)$, a contradiction if $\gamma(\omega)<+\infty$.

\noindent In \cite[Lemma A.1]{sectorialextensions1} we have shown that \eqref{omega7} does imply $\gamma(\omega)=+\infty$.

\item[$(iii)$] Condition \eqref{mixedomega7equ} is clearly stable under \hyperlink{sim}{$\sim$}, which follows by the characterization shown above or can also seen directly since $\omega(t)\ge C^{-1}\sigma(t)-1$ yields $(\omega^{\iota})_{\star}(t^2)=\inf_{s>0}\{\omega(1/s)+t^2s\}\ge C^{-1}\inf_{s>0}\{\sigma(1/s)+t^2Cs\}-1=C^{-1}(\sigma^{\iota})_{\star}(Ct^2)-1\ge C^{-1}(\sigma^{\iota})_{\star}(t^2)-1$ because $(\sigma^{\iota})_{\star}$ is increasing.

\noindent In particular, by \eqref{goodequivalence}, we see that for each $\omega$ as considered in Theorem \ref{mixedomega7} the matrix $\Omega$ satisfies \eqref{Roumieusquare} and/or \eqref{Beurlingsquare} if and only if \eqref{mixedomega7equ} is satisfied for $\omega_{W^{(\lambda)}}$ for some/each $\lambda>0$.

\item[$(iv)$] In general between \eqref{omega7} and \eqref{mixedomega7equ} there is a big difference. As pointed out before, the first condition yields strong non-quasianalyticity for $\omega$, whereas the second one can even be satisfied by (large) quasianalytic weight functions: For this consider the power weights $\omega(t):=t^{\alpha}$, $\alpha\ge 1$, then a straightforward computation yields $$(\omega^{\iota})_{\star}(t)=\left(\alpha^{1/(\alpha+1)}+\frac{1}{\alpha^{\alpha/(\alpha+1)}}\right)t^{\alpha/(\alpha+1)}$$
    and so \eqref{mixedomega7equ} holds true (since $2\alpha/(\alpha+1)\le\alpha\Leftrightarrow 2\le\alpha+1$).
\end{itemize}

So far we have started with a weight function satisfying some standard growth properties, in the next result we will start with a weight sequence $M$ and are interested in the case $\omega\equiv\omega_M$. Recall that for given $M\in\hyperlink{LCset}{\mathcal{LC}}$ the associated weight function $\omega_M$ does have \hyperlink{om0}{$(\omega_0)$}, \hyperlink{om3}{$(\omega_3)$} and \hyperlink{om4}{$(\omega_4)$} (e.g. see \cite[Chapitre I]{mandelbrojtbook}, \cite[Definition 3.1]{Komatsu73} and also \cite{BonetMeiseMelikhov07}).

\begin{proposition}\label{omegaMmixedomega7}
Let $M=(j!m_j)_{j\in\NN}\in\hyperlink{LCset}{\mathcal{LC}}$ and let $\omega_M$ be the associated weight function, $\Omega:=\{W^{(\lambda)}: \lambda>0\}$ shall denote the matrix associated with $\omega_M$. Then $\omega_M$ satisfies \eqref{mixedomega7equ} if and only if
\begin{equation}\label{omegaMmixedomega7equ}
\exists\,C\in\NN_{\ge 1}\,\exists\,D,h>0\text{ such that
}\forall\,j\in\NN\, , \,\,(m_j)^{2C}\le Dh^jm_{Cj}.
\end{equation}
\end{proposition}

Note that \eqref{omegaMmixedomega7equ} is clearly stable under relation \hyperlink{approx}{$\approx$}.

\demo{Proof}
Let $\omega_M$ satisfy \eqref{mixedomega7equ} and w.l.o.g. we can assume $C\in\NN_{\ge 1}$. We follow the ideas in the proof of \cite[Lemma 3.4 $(i)$]{sectorialextensions} (for $M$ instead of $m$). First, for all $j\in\NN$, we get
\begin{align*}
M_{2j}&=\sup_{t>0}\frac{t^{2j}}{\exp(\omega_M(t))}=\frac{e}{e}\sup_{t>0}\frac{(Ht)^{2j}}{\exp(\omega_M(Ht))}\le eH^{2j}\sup_{t>0}\frac{t^{2j}}{\exp(\frac{1}{C}(\omega_M^{\iota})_{\star}(t^2))}
\\&
=eH^{2j}\sup_{t>0}\frac{t^j}{\exp(\frac{1}{C}(\omega_M^{\iota})_{\star}(t))}.
\end{align*}
The supremum in the last expression yields
\begin{align*}
\exp\left(\sup_{t>0}\left\{j\log(t)-\frac{1}{C}(\omega_M^{\iota})_{\star}(t))\right\}\right)&=\exp\left(\sup_{t>0}\left\{j\log(t)-\frac{1}{C}\inf_{s>0}\{\omega_M(1/s)+st\}\right\}\right)
\\&
=\exp\left(\sup_{s,t>0}\left\{j\log(t)-\frac{1}{C}\omega_M(1/s)-\frac{st}{C}\right\}\right).
\end{align*}
By studying for every $j\in\NN$ and $s>0$ fixed the function
$f_{j,s}(t):=j\log(t)-\frac{1}{C}\omega_M(1/s)-\frac{st}{C}$, $t>0$,
one gets that its supremum is given by
$\log\left(\frac{(jC)^j}{(es)^j}\right)-\frac{1}{C}\omega_M(1/s)$ (if $j=0$ we use the convention $0^0:=1$).
Using this we can continue the above estimation for all $j\in\NN$ as follows:
\begin{align*}
eH^{2j}\sup_{t>0}\frac{t^j}{\exp(\frac{1}{C}(\omega_M^{\iota})_{\star}(t))}
&=eH^{2j}\exp\left(\sup_{s>0}\left\{\log\left(\frac{(jC)^j}{(es)^j}\right)-\frac{1}{C}\omega_M(1/s)\right\}\right)
\\
&
  =eH^{2j}\frac{C^jj^j}{e^j}\sup_{s>0}\frac{1}{s^j\exp(\frac{1}{C}\omega_M(1/s))}\\
&=e\left(\frac{CH^2}{e}\right)^jj^j\sup_{s>0}\frac{s^j}{\exp(\frac{1}{C}\omega_M(s))}\\
& =e\left(\frac{CH^2}{e}\right)^jj^j(M_{Cj})^{1/C}.
\end{align*}
Summarizing everything we have shown so far that there exist some
$C_1, h_1>0$ such that for all $j\in\NN$ we get $(M_j)^2\le M_{2j}\le
C_1h_1^jj!(M_{Cj})^{1/C}$ (using for the first estimate that the
log-convexity for $M$ implies that $(M_{j}^{1/j})_{j}$ is increasing) and so
$(M_j)^{2C}\le C_2h_2^jj!^CM_{Cj}$ which is equivalent to
$$ j!^C(m_j)^{2C}\le C_2h_2^j(Cj)!m_{Cj}$$
for all $j \in \NN$. Since by Stirling's formula $(Cj)!$ is growing
like $j!^C$ up to a factor with exponential growth,
 we obtain $(m_j)^{2C}\le C_3h_3^jm_{Cj}$ for all $j \in \NN$ and for
 some constants $C_3$, $h_3$ not depending on $j$, thus \eqref{omegaMmixedomega7equ} is verified.

\medskip

Conversely, assume that \eqref{omegaMmixedomega7equ} is valid. By
going back in the equivalences above, we get $(M_j)^2\le
D_1h_1^jj!(M_{Cj})^{1/C}$ for all $j \in \NN$. If
$\Omega:=\{W^{(\lambda)}: \lambda>0\}$ denotes the matrix associated
with $\omega_M$, then it is known and straightforward to verify that $M\equiv W^{(1)}$ (e.g. see the proof of \cite[Thm. 6.4]{testfunctioncharacterization}) and moreover $W^{(\lambda)}_j=\exp(\frac{1}{\lambda}\varphi^{*}_{\omega}(\lambda j))=(W^{(1)}_{\lambda j})^{1/\lambda}$ for all $\lambda\in\NN_{\ge 1}$. Thus we obtain $(W^{(1)}_j)^2=(M_j)^2\le D_1h_1^jj!(M_{Cj})^{1/C}=D_1h_1^j\widehat{W}^{(C)}_j$ for all $j\in\NN$. Then follow the first part in the proof of Theorem \ref{mixedomega7} with $\lambda=1$ and $\kappa=C$ in order to conclude.
\qed\enddemo


By combining now Proposition \ref{pointwisecharacterizationmatrix}, \eqref{goodequivalence}, Theorem \ref{mixedomega7}, $(iii)$ in the previous observations and Proposition \ref{omegaMmixedomega7} we get the following result.

\begin{corollary}\label{mixedomega7cor1}
Let $\omega$ be given satisfying \hyperlink{om0}{$(\omega_0)$}, \hyperlink{om3}{$(\omega_3)$} and \hyperlink{om4}{$(\omega_4)$}, let $\Omega=\{W^{(\lambda)}: \lambda>0\}$ be the matrix associated with $\omega$. Then the following are equivalent:
\begin{itemize}
\item[$(i)$] $\mathcal{F}_{\{\Omega\}}$ and/or $\mathcal{F}_{(\Omega)}$ is stable under the pointwise product $\odot$,

\item[$(ii)$] $\Omega$ satisfies \eqref{Roumieusquare} and/or \eqref{Beurlingsquare},

\item[$(iii)$] $\omega$ satisfies \eqref{mixedomega7equ},

\item[$(iv)$] some/each $\omega_{W^{(\lambda)}}$ satisfies \eqref{mixedomega7equ},

\item[$(v)$] some/each $W^{(\lambda)}$ satisfies \eqref{omegaMmixedomega7equ}.
\end{itemize}
If $\omega$ has in addition \hyperlink{om1}{$(\omega_1)$}, then
in $(i)$ we have
$\mathcal{F}_{[\omega]}=\mathcal{F}_{[\Omega]}$ and so stability of $\mathcal{F}_{[\omega]}$ under $\odot$ is characterized.
\end{corollary}

On the other hand, starting with a weight sequence satisfying an additional assumption, we have the following characterization.

\begin{corollary}\label{mixedomega7cor2}
Let $M\in\hyperlink{LCset}{\mathcal{LC}}$ be given and satisfying \hyperlink{mg}{$(\on{mg})$}, then the following are equivalent:

\begin{itemize}
\item[$(i)$] $M$ satisfies \eqref{pointwise} (i.e.  $\mathcal{E}_{\{M\}}\subseteq\mathcal{C}^{\omega}$ and/or $\mathcal{E}_{(M)}\subseteq\mathcal{C}^{\omega}$),

\item[$(ii)$] $M$ satisfies \eqref{omegaMmixedomega7equ},

\item[$(iii)$] $\omega_M$ satisfies \eqref{mixedomega7equ},

\item[$(iv)$] $\mathcal{F}_{\{M\}}=\mathcal{F}_{\{\Omega\}}$ and/or $\mathcal{F}_{(M)}=\mathcal{F}_{(\Omega)}$ is stable under the pointwise product.
\end{itemize}
\end{corollary}

\demo{Proof}
Under the assumptions on $M$ we have $\mathcal{F}_{\{M\}}=\mathcal{F}_{\{\Omega\}}$ and/or $\mathcal{F}_{(M)}=\mathcal{F}_{(\Omega)}$ which follows analogously as for the corresponding ultradifferentiable function classes by having the same seminorms, see the proofs of \cite[Cor. 5.8 $(2)$, Lemma 5.9 $(5.11)$]{compositionpaper}. In fact all $W^{(\lambda)}$ are equivalent to $W^{(1)}\equiv M$. Consequently, by combining Proposition \ref{pointwisecharacterizationmatrix} applied to $M\equiv\mathcal{M}$, Theorem \ref{mixedomega7} and finally Proposition \ref{omegaMmixedomega7} we are done.\vspace{6pt}

The equivalence $(i)\Leftrightarrow(ii)$ can also be seen directly as
follows: On the one hand, $(i)\Rightarrow(ii)$ holds by having
$(m_j)^2\le Ab^jm_j$ and so take $C=1$ in
\eqref{omegaMmixedomega7equ}. Conversely, by
assumption $M$ has \hyperlink{mg}{$(\on{mg})$}, i.e. $M_{j+k}\le A_{0}^{j+k}M_jM_k$ for all
$j,k\in\NN$ and some constant $A_{0}$. Consequently,  $m_{j+k}\le
A_{1}^{j+k}m_jm_k$ for all $j,k\in\NN$ and some constant $A_{1}$. By
\eqref{omegaMmixedomega7equ}, we have $(m_j)^{2C}\le Dh^jm_{Cj}$ and
by iteration of $m_{j+k}\le
A_{1}^{j+k}m_jm_k$ , we get $Dh^jm_{Cj}\le B b^{C^2j}(m_j)^C$ and so $(m_j)^2\le
B^{1/C}b^{Cj}m_j$ for some constants $b,B>0$ which is precisely \eqref{pointwise}.

\qed\enddemo

The next result establishes a connection between \eqref{omegaMmixedomega7equ} and the non-quasianalyticity of a sequence $M$.

\begin{lemma}\label{destroyer}
Let $M\in\hyperlink{LCset}{\mathcal{LC}}$ be given such that $\sup_{j\in\NN_{>0}}(m_j)^{1/j}=+\infty$ and \eqref{omegaMmixedomega7equ} holds true. Then $M$ is non-quasianalytic.
\end{lemma}

\demo{Proof}
First, $M$ has \eqref{omegaMmixedomega7equ} if and only if there exist
$C \in \NN_{\geq 1}$ and $D,h>0$ such that $(m_j)^{2C}\le Dh^jm_{Cj}$,
which is equivalent to
$$((M_j)^{1/j})^2\le\frac{j!^{2/j}}{(Cj)!^{1/(Cj)}}D^{1/(Cj)}h^{1/C}(M_{Cj})^{1/(Cj)}.$$
By Stirling's formula $\frac{j!^{2/j}}{(Cj)!^{1/(Cj)}}$ is asymptotically growing like $j\mapsto D_1j$ and so $M$ has \eqref{omegaMmixedomega7equ} if and only if
\begin{equation}\label{destroyerequ}
\exists\,C\in\NN_{\ge 1}\,\exists\,C_1\ge 1\text{ such that
}\forall\,j\in\NN_{>0}\, , \,\,((M_j)^{1/j})^2\le C_1j(M_{Cj})^{1/(Cj)}.
\end{equation}
Note that the assumption $\sup_{j\in\NN_{>0}}(m_j)^{1/j}=+\infty$
implies that in \eqref{destroyerequ} we have $C\ge 2$: indeed, the case $C=1$ would yield \eqref{pointwise} and so $\sup_{j\in\NN_{>0}}(m_j)^{1/j}<+\infty$, hence a contradiction.

Since we have $\sup_{j\in\NN_{>0}}(m_j)^{1/j}=+\infty$, for all $A\ge 1$ there does exist a number $q_A\in\NN_{\ge 1}$ (which can be chosen minimal) such that we get
$(m_{q_A})^{1/q_A}\ge A$, or equivalently $(M_{q_A})^{1/q_A}\ge A(q_A!)^{1/q_A}$. Thus, by a consequence of Stirling's formula, we obtain $(M_{q_A})^{1/q_A}\ge\frac{Aq_A}{e}$ and so also $\frac{eCC_1}{A}\ge\frac{CC_1q_A}{(M_{q_A})^{1/q_A}}$ follows with $C$ and $C_1$ denoting the constants arising in \eqref{destroyerequ} (which are not depending on given $q_A$).

Let now $A\ge 1$ be chosen sufficiently large in order to have $\frac{eCC_1}{A}<1$ and set $q:=q_A$. By the above we see that $\frac{CC_1q}{(M_q)^{1/q}}<1$ holds true.

Since $M\in\hyperlink{LCset}{\mathcal{LC}}$ we have that
$j\mapsto(M_j)^{1/j}$ is increasing. As we will see this property is sufficient to conclude and for convenience we put now $L_j:=(M_j)^{1/j}$. For the sum under consideration we estimate by
$$\sum_{j\ge q}\frac{1}{(M_j)^{1/j}}=\sum_{k=0}^{+\infty}\sum_{j=C^kq}^{C^{k+1}q-1}\frac{1}{L_j}\le\sum_{k=0}^{+\infty}\frac{C^{k+1}q-C^kq}{L_{C^kq}}=\sum_{k=0}^{+\infty}\frac{C^kq(C-1)}{L_{C^kq}}.$$
Now, by iterating \eqref{destroyerequ}, we have for every $k\in\NN_{\ge 1}$
\begin{align*}
\frac{C^kq(C-1)}{L_{C^kq}}&
 \le  qC^{k+1}\frac{1}{L_{CC^{k-1}q}}\\
&\le  qC^{k+1}C_1C^{k-1}q\frac{1}{(L_{C^{k-1}q})^2}\\
&= q^2C^{2^1(k-0)}C_1^2\frac{1}{(L_{CC^{k-2}q})^2}
\\&
\le
    q^2C^{2^1(k-0)}C_1^2C_1^2C^{2k-4}q^2\frac{1}{(L_{C^{k-2}q})^4}\\
& = q^4C^{2^2(k-1)}C_1^4\frac{1}{(L_{CC^{k-3}q})^4}
\\&
\le q^4C^{2^2(k-1)}C_1^4
    C_1^4(C^{k-3})^4q^4\frac{1}{(L_{C^{k-3}q})^8}\\
& =q^8C_1^8C^{2^3(k-2)}\frac{1}{(L_{C^{k-3}q})^8}
\\&
\le\dots\le\frac{q^{2^k}C_1^{2^k}C^{2^k}}{(L_{q})^{2^k}},
\end{align*}
where we have used that for all natural numbers $i,k$ with $1\le i\le k-1$ we get $(C^{k-(i+1)}q)^{2^{i}}=q^{2^{i}}C^{2^{i}(k-(i+1))}$ and $C^{2^i(k-(i-1))}C^{2^{i}(k-(i+1))}=C^{2^{i+1}(k-i)}$.

Finally, if $k=0$, then $\frac{C^kq(C-1)}{L_{C^kq}}\le\frac{qC}{L_q}\le\frac{qCC_1}{L_q}$ and gathering everything we have shown now
$$\sum_{k=0}^{+\infty}\frac{C^kq(C-1)}{L_{C^kq}}\le\sum_{k=0}^{+
\infty}\Bigl(\underbrace{\frac{qCC_1}{L_q}}_{<1}\Bigr)^{2^k}<+\infty,$$
which proves the non-quasianalyticity for $M$ as desired.
\qed\enddemo

Using the above Lemma we can prove now the final statement of this section showing that the problem of algebrability cannot be considered within the quasianalytic weight function setting.

\begin{theorem}\label{destroyercor}
Let $\omega$ satisfying \hyperlink{om0}{$(\omega_0)$}, \hyperlink{om2}{$(\omega_2)$}, \hyperlink{om3}{$(\omega_3)$}, \hyperlink{om4}{$(\omega_4)$} and $\liminf_{t\rightarrow+ \infty}\frac{\omega(t)}{t}=0$ be given. Assume that $\omega$ has in addition the characterizing condition \eqref{mixedomega7equ} (resp. equivalently $\mathcal{F}_{\{\Omega\}}$ and/or $\mathcal{F}_{(\Omega)}$ is stable under the pointwise product $\odot$), then $\omega$ has to be non-quasianalytic, i.e. condition \eqref{omegaQ} is violated.
\end{theorem}

\demo{Proof}
Let $\Omega=\{W^{(\lambda)}: \lambda>0\}$ be the matrix associated
with $\omega$. We apply Lemma \ref{destroyer} to some/each sequence
$W^{(\lambda)}$ which can be done by the assumptions on $\omega$
and the equivalences obtained in Corollary \ref{mixedomega7cor1} above. Then $W^{(\lambda)}$ has \hyperlink{mnq}{$(\on{nq})$} and so $\omega$ does not enjoy \eqref{omegaQ} (recall that this last step holds by \cite[Lemma 4.1]{Komatsu73} and \eqref{goodequivalence}).
\qed\enddemo

Note that this result deals with a property of the associated matrix $\Omega$ and \hyperlink{om1}{$(\omega_1)$} is not required necessarily. If $\omega$ has in addition \hyperlink{om1}{$(\omega_1)$}, then we have $\mathcal{F}_{[\omega]}=\mathcal{F}_{[\Omega]}$ in Theorem \ref{destroyercor}.

\medskip

\textbf{Acknowledgements.}
\medskip
The authors wish to thank the referee for his comments which have improved the presentation and the structure of this work. The authors also wish to thank Javier Jim\'enez-Garrido and Javier Sanz from the Universidad de Valladolid for their helping discussions concerning the results and proofs of Section \ref{stabilitypointwiseweightfunction}.

\bibliographystyle{plain}
\bibliography{Bibliography}
\end{document}